\documentclass[11pt]{article}

\usepackage{amsmath,amssymb,amsthm,amsfonts,mathabx,enumitem}

\def\newcaption#1{\vspace{-.1in}\caption{#1}}
\def\tbl{\begin{table}[h]\begin{small}\begin{center}}
\def\tblH{\begin{table}[H]\begin{small}\begin{center}}
\def\endtbl#1{\end{center}\end{small}\label{#1}\vspace{-.2in}\end{table}}

\newtheorem{thm}{Theorem}
\theoremstyle{remark}
\newtheorem{rmk}{Remark}

\def\BE#1{\begin{equation}\label{#1}}
\def\EE{\end{equation}}
\def\eref#1{(\ref{#1})}

\def\lr#1{\langle{#1}\rangle}
\def\blr#1{\big\langle{#1}\big\rangle}
\def\ov#1{\overline{#1}}

\def\sf#1{\textsf{#1}}
\def\wch#1{\widecheck{#1}}
\def\tn#1{\textnormal{#1}} 

\def\wt#1{\widetilde{#1}}

\def\lra{\longrightarrow}

\def\al{\alpha}

\def\de{\delta}
\def\ep{\epsilon}
\def\io{\iota}

\def\om{\omega}

\def\vph{\varphi}

\def\Si{\Sigma}

\def\C{\mathbb C}
\def\fd{\mathfrak d}
\def\bF{\mathbb F}
\def\cJ{\mathcal J}

\def\cN{\mathcal N}
\def\fo{\mathfrak o}

\def\os{\mathfrak{os}}
\def\P{\mathbb P}
\def\Q{\mathbb Q}
\def\R{\mathbb R}
\def\fs{\mathfrak s}

\def\Z{\mathbb Z}
\def\cH{\mathcal H}
\def\cP{\mathcal P}

\def\a{\mathbf a}
\def\b{\mathbf b}

\def\m{\mathbf m}

\def\nd{\tn{d}}
\def\Im{\tn{Im}}
\def\pt{\tn{pt}}
\def\PD{\tn{PD}}

\def\part{\partial}

\def\1{\mathbf 1}
\def\2{\mathbf 2}
\def\3{\mathbf 3}
\def\4{\mathbf 4}

\def\prt{\partial}
\def\st{\bigstar}

\begin{document}

\title{WDVV-Type Relations for Welschinger's Invariants: Applications}
\author{Xujia Chen and Aleksey Zinger\thanks{Supported by NSF grants DMS 1500875 and 1901979
and Simons collaboration grant 587036}}
\date{\today}

\maketitle

\begin{abstract}
\noindent
We first recall Solomon's relations for Welschinger's invariants
counting real curves in real symplectic fourfolds,
announced in~\cite{Jake2} and established in~\cite{RealWDVV},
and the WDVV-style relations for Welschinger's invariants
counting real curves in real symplectic sixfolds with some symmetry
established in~\cite{RealWDVV3}.
We then explicitly demonstrate that in some important cases
(projective spaces with standard conjugations, 
real blowups of the projective plane,
and two- and three-fold products of the one-dimensional projective space
with two involutions each), 
these relations provide complete recursions determining all Welschinger's invariants
from basic input.
We include extensive tables of Welschinger's invariants
in low degrees obtained from these recursions with {\it Mathematica}.
These invariants provide lower bounds for counts of real rational curves,
including with curve insertions in smooth algebraic threefolds.
\end{abstract}

\tableofcontents

\section{Introduction}
\label{intro_sec}

Let $(X,\om,\phi)$ be a compact connected real symplectic manifold of (real) dimension~$2n$.
The fixed locus~$X^{\phi}$ of the anti-symplectic involution~$\phi$ on~$X$
is then a Lagrangian submanifold of~$(X,\om)$.
Let 
\begin{gather*}
H_2^{\phi}(X)=\big\{B\!\in\!H_2(X;\Z)\!:\phi_*B\!=\!-B\big\},\qquad
H^*(X)^{\phi}_{\pm}=\big\{\mu\!\in\!H^*(X)\!:\phi^*\mu\!=\!\pm\mu\big\}, \\
\fd\!:H_2(X;\Z)\lra H_2^{\phi}(X), ~~ \fd(B)=B\!-\!\phi_*(B).
\end{gather*}
We denote by $\cJ_{\om}$ the space of $\om$-compatible (or -tamed) 
almost complex structures~$J$ on~$X$ and by
$\cJ_{\om}^{\phi}\!\subset\!\cJ_{\om}$ the subspace of almost complex structures~$J$
such that \hbox{$\phi^*J\!=\!-J$}.
Let 
$$c_1(X,\om)\equiv c_1(TX,J)\in H^2(X)$$
be the first Chern class of $TX$ with respect to some $J\!\in\!\cJ_{\om}$;
it is independent of such a choice.
Define
$$\ell_{\om}\!:H_2(X;\Z)\lra\Z, \quad 
\ell_{\om}(B)=\blr{c_1(X,\om),B}\!+\!n\!-\!3.$$
The paradigmatic example of a real symplectic manifold is the complex projective space~$\P^n$
with the Fubini-Study symplectic form and the standard conjugation
$$\tau_n\!:\P^n\lra\P^n, \qquad 
\tau_n\big([Z_0,Z_1,\ldots,Z_n]\big)=\big[\ov{Z_0},\ov{Z_1},\ldots,\ov{Z_n}\big].$$

\vspace{.2in}

For $J\!\in\!\cJ_{\om}$ and $B\!\in\!H_2(X;\Z)$, 
a subset $C\!\subset\!X$ is a \sf{genus~0} (or \sf{rational}) 
\sf{irreducible $J$-holomorphic degree~$B$ curve}
if there exists a simple (not multiply covered) \hbox{$J$-holomorphic} map 
\BE{Cudfn_e}u\!:\P^1\lra X \qquad\hbox{s.t.}\quad  C=u(\P^1),~~u_*[\P^1]=B.\EE
If $J\!\in\!\cJ_{\om}^{\phi}$, such a curve~$C$ is \sf{real} if 
in addition \hbox{$\phi(C)\!=\!C$};
if so, then $B\!\in\!H_2^{\phi}(X;\Z)$.\\

Invariant signed counts of real rational $J$-holomorphic curves in compact 
real symplectic fourfolds and sixfolds, now known as \sf{Welschinger's invariants},
were defined in~\cite{Wel4,Wel6} and interpreted in terms of moduli spaces of $J$-holomorphic maps
from the two-disk in~\cite{Jake}. 
An adaptation of the interpretation of~\cite{Jake} in terms of real $J$-holomorphic maps
from the Riemann sphere~$\P^1$ later appeared in~\cite{Penka2} and
was reformulated in terms of degrees of relatively oriented pseudocycles
in \cite{RealWDVV,RealWDVV3}.
The moduli-theoretic perspectives on Welschinger's invariants
lead to the WDVV-type relations for them established in~\cite{RealWDVV,RealWDVV3}.\\

We recall the relevant versions of definitions of Welschinger's invariants of
real symplectic fourfolds and sixfolds in Sections~\ref{genWel_sec} and~\ref{genWel3_sec}, 
respectively.
We conclude both sections with statements of WDVV-type relations for
Welschinger's invariants in the respective settings;
see Theorems~\ref{WelWDVV_thm} and~\ref{WelWDVV3_thm}, respectively.
In the case of~$(\P^2,\tau_2)$,
the two relations of Theorem~\ref{WelWDVV_thm} restrict to the two relations of~\cite[p13]{Jake2},
which determine all Welschinger's invariants of~$(\P^2,\tau_2)$ from basic input; 
see Section~\ref{P2_sec}.
We explicitly demonstrate in Section~\ref{blowups_sec} that 
the first relation of Theorem~\ref{WelWDVV_thm} and elementary geometric consideration
reduce all Welschinger's invariants of real blowups $(\P^2_{r,s},\tau_{r,s})$ of~$(\P^2,\tau_2)$ 
to the standard Gromov-Witten invariants of
blowups~$\P^2_k$ of~$\P^2$ and Welschinger's invariants of~$(\P^2,\tau_2)$.
In Section~\ref{P1P1_sec}, we apply Theorem~\ref{WelWDVV_thm}
to $\P^1\!\times\!\P^1$ with two natural orientation-reversing involutions
and show that it determines all Welschinger's invariants from basic input 
in theses cases as well.\\

In Sections~\ref{P3_sec} and~\ref{P1P1P1_sec}, 
we explicitly show that the relations of Section~\ref{genWel3_sec}
similarly determine all Welschinger's invariants of~$(\P^3,\tau_3)$
and $(\P^1)^3$ with two distinct conjugations from some basic inputs.
However, new nuances arise in these settings. 
First, Welschinger's invariants of a real symplectic sixfold $(X,\om,\phi)$
depend on the choice of a Spin-structure on~$X^{\phi}$
in all perspectives on these invariants.
Second, the WDVV-type relations of Section~\ref{genWel3_sec} for these invariants
involve signed counts of real curves through constraints of dimension~2.
Such counts are not part of Welschinger's original definition in~\cite{Wel6};
they arise from the moduli-theoretic perspective on Welschinger's invariants 
introduced in~\cite{Jake}
and the use of a symmetry of $(X,\om,\phi)$ introduced in~\cite{RealWDVV3}.
Finally, the WDVV-type relations of Section~\ref{genWel3_sec} determine
Welschinger's invariants of $(\P^3,\tau_3)$ from a single input,
as shown in \cite[Section~4.1.4]{Adam};
the sign of this input depends on the choice of a Spin-structure on
the fixed locus $\R\P^3\!\subset\!\P^3$ of~$\tau_3$.
On the other hand, the determination of Welschinger's invariants of $(\P^1)^3$
with either of the two natural conjugations~$\phi$ requires a basic input beyond those
that can be directly changed by changing the Spin structure on the fixed locus of~$\phi$.
We determine this input via a reduction to~$(\P^1)^2$.
The tables at the end provide low-degree Welschinger's invariants based 
on the formulas in Sections~\ref{P2_sec}-\ref{blowups_sec}, \ref{P3_sec}, and~\ref{P1P1P1_sec};
the {\it Mathematica programs} implementing these formulas are
available from the Wolfram Foundation's {\it Notebook Archive}.\\

The authors would like to thank E.~Brugall\'e for suggestions on
a predecessor to the present paper.

\section{Solomon's relations for real symplectic fourfolds}
\label{genWel_sec}

Let $(X,\om,\phi)$ be a compact real symplectic fourfold with connected fixed locus~$X^{\phi}$.
For \hbox{$B\!\in\!H_2(X;\Z)$}, define
$$\lr{B}_l=\begin{cases}1,&\hbox{if}~2l\!=\!\ell_{\om}(B)\!-\!1;\\
0,&\hbox{otherwise}.\end{cases}$$
Suppose $J\!\in\!\cJ_{\om}^{\phi}$, $C\!\subset\!X$ is a real curve,
and $u$ is as in~\eref{Cudfn_e}.
A point $x\!\in\!C$ is a \sf{simple node} if
$$\big|u^{-1}(x)\big|= 2 \qquad\hbox{and}\qquad
\bigoplus_{z\in u^{-1}(x)}\!\!\!\!\!\Im\,\nd_zu=T_xX\,.$$
We denote by $\de_E(C)$ the number of nodes of~$C$ that are isolated points
of~$C\!\cap\!X^{\phi}$.\\

Suppose $B\!\in\!H_2(X;\Z)$ and $l\!\in\!\Z^{\ge0}$ are such~that
\BE{dimcond_e} k\equiv \ell_{\om}(B)\!-\!2l\in\Z^{\ge0}\,.\EE
For a generic $J\!\in\!\cJ_{\om}^{\phi}$, 
there are then only finitely many rational irreducible real $J$-holomorphic 
degree~$B$ curves $C\!\subset\!X$ passing through  $k$ points in~$\wch{X}^{\phi}$
and $l$ points in $X\!-\!X^{\phi}$ in general position. 
According to \cite[Thm.~0.1]{Wel4}, the~sum
$$N_{B,l}^{\phi}\equiv \sum_C{(-1)}^{\de_E(C)}$$
over the set of these curves
is independent of the choices of~$J$ and the points.
If the number~$k$ in~\eref{dimcond_e} is negative, 
we set $N_{B,l}^{\phi}\!=\!0$.\\

If $B\!\in\!H_2(X;\Z)$,  $\ell_{\om}(B)\!\ge\!0$, and $J\!\in\!\cJ_{\om}$, 
there are only finitely many rational irreducible $J$-holomorphic 
degree~$B$ curves~$C$ passing through $\ell_{\om}(B)$ points in general position. 
The number of such curves counted with appropriate signs is independent of the choices of~$J$ 
and the points. 
This is the standard (complex) genus~0 degree~$B$ \sf{Gromov-Witten invariant} of~$(X,\om)$
with $\ell_{\om}(B)$ point insertions;
we denote it by~$N_B^X$.
If $\ell_{\om}(B)\!<\!0$, we set $N_B^X\!=\!0$.\\

For $B,B'\!\in\!H_2(X;\Z)$, we denote by $B\!\cdot_X\!\!B'\!\in\!\Z$  
the homology intersection product of~$B$ with~$B'$
and by $B^2\!\in\!\Z$ the self-intersection number of~$B$. 
Since $\phi$ is an anti-symplectic involution, 
$\fd(B)^2\!\in\!2\Z$ for every $B\!\in\!H_2(X;\Z)$.
The theorem below follows immediately from the proof of \cite[Thm.~1.1]{RealWDVV},
which establishes similar relations for signed counts of real $J$-holomorphic maps
as in \cite{Jake,Penka2}, and \cite[Thm~12.1]{SpinPin},
which compares the signs of these counts with the signs in Welschinger's definition in~\cite{Wel4}.

\begin{thm}\label{WelWDVV_thm}
Suppose $(X,\om,\phi)$ is a compact real symplectic fourfold, 
$H_1,H_2\!\in\!H^2(X)^{\phi}_-$, $B\!\in\!H_2(X)$, and $l\!\in\!\Z^{\ge0}$. 
\begin{enumerate}[label=(SWDVV\arabic*),leftmargin=*]

\item\label{Wel12rec_it} If $l\!\ge\!1$ and $\ell_{\om}(B)\!-\!2l\!\ge\!1$ 
(i.e.~$k\!\ge\!1$), then
\begin{equation*}\begin{split}
&\hspace{-.8in}
\lr{H_1H_2,X}N_{B,l}^{\phi}
= (-1)^{l+B^2/2}2^{l-3}\lr{B}_l\lr{H_1,B}\lr{H_2,B}
\hspace{-.15in}\sum_{\begin{subarray}{c}B'\in H_2(X)\\ \fd(B')=B\end{subarray}}
\hspace{-.18in}N_{B'}^X \\
&\hspace{-.7in}
+\!\!\!\!\!\!\sum_{\begin{subarray}{c}B_0,B'\in H_2(X)-\{0\}\\ B_0+\fd(B')=B \end{subarray}}
\hspace{-.4in}(-1)^{l_\om(B')+\fd(B')^2/2}2^{l_\om(B')}\!\big(B_0\!\!\cdot_X\!\!B'\big)
\lr{H_1,B'}\lr{H_2,B'}\binom{l\!-\!1}{\ell_{\om}(B')}
N_{B'}^XN_{B_0,l-1-\ell_{\om}(B')}^{\phi}\\
&\hspace{-.7in}+\!\!\!\!\!\!\!
\sum_{\begin{subarray}{c}B_1,B_2\in H_2(X)-\{0\}\\ B_1+B_2=B \\ 
l_1+l_2=l-1,\,l_1,l_2\ge0\end{subarray}} \hspace{-.4in}
\lr{H_1,B_1}\binom{l\!-\!1}{l_1}\!\!
\Bigg(\!\lr{H_2,B_2}\binom{\ell_{\om}(B)\!-\!2l\!-\!1}{\ell_{\om}(B_1)\!-\!2l_1\!-\!1}
-\lr{H_2,B_1}\binom{\ell_{\om}(B)\!-\!2l\!-\!1}{\ell_{\om}(B_1)\!-\!2l_1}\!\!\!\Bigg)\!
N_{B_1,l_1}^{\phi}N_{B_2,l_2}^{\phi}.
\end{split}\end{equation*}

\item\label{Wel03rec_it} If $l\!\ge\!2$, then
\begin{equation*}\begin{split}
&\hspace{-.8in}
\lr{H_1H_2,X}N_{B,l}^{\phi}=-\!\!\!\!\!
\sum_{\begin{subarray}{c}B_0,B'\in H_2(X)-\{0\}\\ B_0+\fd(B')=B\end{subarray}}
\hspace{-.4in} (-1)^{l_\om(B')+\fd(B')^2/2}
2^{l_\om(B')-1}\!\big(B_0\!\!\cdot_X\!\!B'\big)\lr{H_1,B'}\\
&\hspace{.7in}\times\!\Bigg(\!\!\lr{H_2,B_0}\binom{l\!-\!2}{\ell_{\om}(B')\!-\!1}
-2\lr{H_2,B'}\binom{l\!-\!2}{\ell_{\om}(B')}\!\!\Bigg)
N_{B'}^XN_{B_0,l-1-\ell_{\om}(B')}^{\phi}\\
&\hspace{-.7in}+\!\!\!\!\!\!\!
\sum_{\begin{subarray}{c}B_1,B_2\in H_2(X)-\{0\}\\ B_1+B_2=B\\ 
l_1+l_2=l-2,\,l_1,l_2\ge0\end{subarray}} \hspace{-.4in}\lr{H_2,B_1}
\binom{l\!-\!2}{l_1}\!\!
\Bigg(\!\lr{H_1,B_2}\binom{\ell_{\om}(B)\!-\!2l}{\ell_{\om}(B_1)\!-\!2l_1\!-\!1}
-\lr{H_1,B_1}\binom{\ell_{\om}(B)\!-\!2l}{\ell_{\om}(B_1)\!-\!2l_1}\!\!\!\Bigg)\!
N_{B_1,l_1}^{\phi}N_{B_2,l_2+1}^{\phi}.
\end{split}\end{equation*}

\end{enumerate}
\end{thm}

\begin{rmk}\label{WelWDVV_rmk}
Theorem~\ref{WelWDVV_thm} extends to 
compact real symplectic fourfolds $(X,\om,\phi)$ with disconnected fixed loci~$X^{\phi}$
and applies with finer notions of the curve degree~$B$;
see Theorem~1.1 and Remark~1.3 in~\cite{RealWDVV}.
This makes no difference for the examples of Sections~\ref{P2_sec}-\ref{blowups_sec}.
Theorem~1.1 in~\cite{RealWDVV} contains another WDVV-type relation; 
it involves three divisior insertions.
Since it is not needed for the present purposes, we do not state it here.
\end{rmk}

\section{Projective plane~$\P^2$}
\label{P2_sec}

The fixed locus of the standard conjugation 
$$\tau_2\!:\P^2\lra\P^2, \qquad \tau_2\big([Z_0,Z_1,Z_2]\big)
=\big[\ov{Z_0},\ov{Z_1},\ov{Z_2}\big],$$
on the complex projective plane is the real projective plane~$\R\P^2$.
The group \hbox{$H_2^{\tau_2}(\P^2)\!=\!H_2(\P^2;\Z)$}
is identified with~$\Z$ via the standard generator $L\!=\![\P^1]$.\\

For $d\!\in\!\Z^+$ and $l\!\in\!\Z^{\ge0}$, we~set
$$N_d^{\P^2}=N_{dL}^{\P^2}, \qquad N_{d,l}^{\tau_2}=N_{dL,l}^{\tau_2}, \qquad 
\lr{d}_l=\lr{dL}_l\equiv \begin{cases}1,&\hbox{if}~2l\!=\!3d\!-\!2;\\
0,&\hbox{otherwise}.\end{cases}$$
Since there is a unique line~$\P^1$ through every pair of points in~$\P^2$, 
$$N_1^{\P^2},N_{1,0}^{\tau_2},N_{1,1}^{\tau_2}=1.$$
Since there is a unique conic through five general points in~$\P^2$
and this conic is an embedded~$\P^1$,
$$N_2^{\P^2},N_{2,0}^{\tau_2},N_{2,1}^{\tau_2},N_{2,2}^{\tau_2}=1.$$

\vspace{.15in}

Kontsevich's recursion for the standard Gromov-Witten invariants of $\P^2$ can be written~as
$$N_d^{\P^2}=\frac{1}{6(d\!-\!1)}\sum_{\begin{subarray}{c}d_1+d_2=d\\ d_1,d_2\ge1\end{subarray}}
\Bigg(\!d_1d_2\!-2\frac{(d_1\!-\!d_2)^2}{3d-2}\!\Bigg)
\binom{3d\!-\!2}{\!3d_1\!-\!1\!}\! d_1d_2N_{d_1}^{\P^2}N_{d_2}^{\P^2}\,.$$
Taking $H_1,H_2\!\in\!H^2(\P^2)$ to be the Poincare duals of~$L$ in
the two formulas of Theorem~\ref{WelWDVV_thm}, we obtain
\begin{enumerate}[leftmargin=49pt]

\item[{\cite[(10)]{Jake2}}] If $d\!\ge\!2$, $l\!\ge\!1$, and $3d\!-\!2l\!\ge\!2$ 
(i.e.~$k\!\ge\!1$), then
\begin{equation*}\begin{split}
{}\hspace{-.8in}
N_{d,l}^{\tau_2}=-(-2)^{3d/2-4}\lr{d}_ld^2N_{d/2}^{\P^2}
+\!\!\!\sum_{\begin{subarray}{c}d_0+2d'=d\\ d_0,d'\ge1\end{subarray}}
\hspace{-.15in} (-2)^{3d'-1}d_0d'^3\binom{l\!-\!1}{3d'\!-\!1}
N_{d'}^{\P^2}N_{d_0,l-3d'}^{\tau_2}\hspace{.3in}&\\
+\!\!\!\sum_{\begin{subarray}{c}d_1+d_2=d\\ l_1+l_2=l-1\\
d_1,d_2\ge1,\,l_1,l_2\ge0\end{subarray}}  \hspace{-.2in}
\binom{l\!-\!1}{l_1}\!\!\Bigg(\!d_1d_2\binom{3d\!-\!2l\!-\!2}{3d_1\!-\!2l_1\!-\!2}
\!-\!d_1^2\binom{3d\!-\!2l\!-\!2}{3d_1\!-\!2l_1\!-\!1}\!\!\Bigg)\!
N_{d_1,l_1}^{\tau_2}\!N_{d_2,l_2}^{\tau_2}&\,.
\end{split}\end{equation*}

\item[{\cite[(9)]{Jake2}}] If $d\!\ge\!2$ and $l\!\ge\!2$, then
\begin{equation*}\begin{split}
{}\hspace{-.8in}
N_{d,l}^{\tau_2}=\sum_{\begin{subarray}{c}d_0+2d'=d\\ d_0,d'\ge1\end{subarray}}
\hspace{-.15in} (-2)^{3d'-2}d_0d'^2
\Bigg(\!\!d_0\!\binom{l\!-\!2}{3d'\!-\!2}\!-\!2d'\binom{l\!-\!2}{3d'\!-\!1}\!\!\Bigg)
N_{d'}^{\P^2}N_{d_0,l-3d'}^{\tau_2}\hspace{.8in}&\\
+\!\!\!\sum_{\begin{subarray}{c}d_1+d_2=d\\ l_1+l_2=l-2\\
d_1,d_2\ge1,\,l_1,l_2\ge0\end{subarray}}  \hspace{-.2in}
\binom{l\!-\!2}{l_1}\!\Bigg(\!d_1d_2\binom{3d\!-\!2l\!-\!1}{3d_1\!-\!2l_1\!-\!2}
\!-\!d_1^2\binom{3d\!-\!2l\!-\!1}{3d_1\!-\!2l_1\!-\!1}\!\!\Bigg)
\!N_{d_1,l_1}^{\tau_2}\!N_{d_2,l_2+1}^{\tau_2}&\,.
\end{split}\end{equation*}

\end{enumerate}

\vspace{.10in}

The $(d,l)\!=\!(2,2)$ cases of these recursions reduce to
$$N_{2,2}^{\tau_2}=2N_1^{\P^2}\!-\!N_{1,1}^{\tau_2}N_{1,0}^{\tau_2}, \qquad
N_{2,2}^{\tau_2}=N_{1,1}^{\tau_2}N_{1,0}^{\tau_2}\,.$$
Taking the difference between the above recursions for $N_{d+1,2}^{\tau_2}$
with $d\!\ge\!2$ as suggested at the end of~\cite{Jake2}, we obtain
$$N_{d,0}^{\tau_2}=\!\sum_{\begin{subarray}{c}d_1+d_2=d\\ d_1,d_2\ge1\end{subarray}} \!\!
\Bigg(\!d_1\!\binom{3d\!-\!3}{3d_1\!-\!2}\!-\!(d_2\!+\!1)\!\binom{3d\!-\!3}{3d_1\!-\!3}\!\!\Bigg)
N_{d_1,0}^{\tau_2}N_{d_2+1,1}^{\tau_2}  \qquad\forall\,d\!\ge\!2.$$
The three recursions above determine all numbers $N_{d,l}^{\tau_2}$ with $d\!\ge\!2$.
The low-degree numbers obtained from these recursions and shown in Table~\ref{WelP2_tbl}
agree with \cite[Corollary~6]{IKS09} and \cite[Section~7.3]{ABL}.

\section{Quadric surfaces $\P^1\!\times\!\P^1$}
\label{P1P1_sec}

Let $\tau\!:\P^1\!\lra\!\P^1$ be the standard conjugation on~$\P^1$;
its fixed locus is $\R\P^1\!\approx\!S^1$.
The fixed loci of the anti-holomorphic involutions 
$$\tau_{1,1},\tau_{1,1}'\!:\P^1\!\times\!\P^1\lra\P^1\!\times\!\P^1, \quad
\tau_{1,1}(z_1,z_2)=\big(\tau(z_1),\tau(z_2)\big),~~
\tau_{1,1}'(z_1,z_2)=\big(\tau(z_2),\tau(z_1)\big),$$
are the two-torus $\R\P^1\!\times\!\R\P^1$ and the anti-diagonal two-sphere
$$\big(\P^1\!\times\!\P^1\big)^{\tau_{1,1}'}
=\big\{\big(z,\tau_1(z)\!\big)\!:z\!\in\!\P^1\big\}\approx\P^1\,,$$
respectively.
The groups $H_2^{\tau_{1,1}}(\P^1\!\times\!\P^1)\!\approx\!H_2(\P^1\!\times\!\P^1)$ 
are identified with~$\Z\!\oplus\!\Z$
 via the standard generators
$$L_1\equiv\big[\P^1\!\times\!\pt\big] \qquad\hbox{and}\qquad
L_2\equiv[\pt\!\times\!\P^1].$$
The group $H_2^{\tau_{1,1}'}(\P^1\!\times\!\P^1)$ is the diagonal subgroup
$$\big\{d_1L_1\!+\!d_2L_2\!:d_1,d_2\!\in\!\Z^2,\,d_1\!=\!d_2\big\}
\approx\Z\!\subset\!\Z\!\times\!\Z.$$

\vspace{.15in}

For $a,b\!\in\!\Z^{\ge0}$ with $(a,b)\!\neq\!(0,0)$, $d\!\in\!\Z^+$, and $l\!\in\!\Z^{\ge0}$, we~set
\begin{gather*}
N_{a,b}^{\P^1\!\times\P^1}=N_{aL_1+bL_2}^{\P^1\!\times\P^1},\quad
N_d^{\P^1\!\times\P^1}=
\sum_{\begin{subarray}{c}a+b=d\\ a,b\ge0\end{subarray}}\!\!\!\!N_{a,b}^{\P^1\!\times\P^1}\,,
\quad N_{(a,b),l}^{\tau_{1,1}}=N_{aL_1+bL_2,l}^{\tau_{1,1}}, 
\quad N_{d,l}^{\tau_{1,1}'}=N_{dL_1+dL_2,l}^{\tau_{1,1}'}\,,\\
\lr{a,b}_l= \begin{cases}1,&\hbox{if}~a,b\!\in\!2\Z,\,l\!=\!a\!+\!b\!-\!1;\\
0,&\hbox{otherwise};\end{cases} \quad
\lr{d}_l= \begin{cases}1,&\hbox{if}~l\!=\!2d\!-\!1;\\
0,&\hbox{otherwise}.\end{cases}
\end{gather*}
By symmetry, 
$$N_{a,b}^{\P^1\!\times\P^1}=N_{b,a}^{\P^1\!\times\P^1}, \qquad
N_{(a,b),l}^{\tau_{1,1}}=N_{(b,a),l}^{\tau_{1,1}}.$$
An irreducible degree $(1,b)$-curve in $\P^1\!\times\!\P^1$ 
corresponds to the graph of a ratio of two degree~$b$ polynomials on~$\C$.
Since every such rational function is determined by its values at $2b\!+\!1$ points,
\BE{P1P1_e2}N_{1,b}^{\P^1\!\times\P^1},
N_{(1,b),0}^{\tau_{1,1}},\ldots,N_{(1,b),b}^{\tau_{1,1}},
N_{1,0}^{\tau_{1,1}'},N_{1,1}^{\tau_{1,1}'}=1 \qquad\forall\,b\!\in\!\Z^{\ge0}.\EE
Since a degree $(a,0)$-curve in $\P^1\!\times\!\P^1$ is a degree~$a$ cover of 
a horizontal section,
$$N_{a,0}^{\P^1\!\times\P^1},N_{(a,0),l}^{\tau_{1,1}}=0 \quad\forall\,a\!\ge\!2.$$

\vspace{.15in}
The complex WDVV recursions for $\P^1\!\times\!\P^1$ can be written as
\begin{gather*}
N_{0,0}\equiv0,\qquad N_{1,0},N_{0,1}=1, \\
N_{a,b}=\frac12\!\!\!\!\!\! \sum_{\begin{subarray}{c}a_1+a_2=a\\ 
b_1+b_2=b\\ a_1,a_2,b_1,b_2\ge0 \end{subarray}} \hspace{-.26in}
(a_1b_2\!+\!a_2b_1)
(a_1\!+\!b_1)\!\bigg(\!\!(a_2\!+\!b_2)\!\binom{2a\!+\!2b\!-\!4}{2a_1\!+\!2b_1\!-\!2}
\!-\!(a_1\!+\!b_1)\!\binom{2a\!+\!2b\!-\!4}{2a_1\!+\!2b_1\!-\!1}\!\!\bigg)
N_{a_1,b_1}N_{a_2,b_2}\,.
\end{gather*}

\subsection{The twisted involution}
\label{P1P1b_subs}

Taking $H_1,H_2\!\in\!H^2(\P^1\!\times\!\P^1)$ to be the Poincare duals 
of $(L_1\!+\!L_2)/\sqrt2$ in
the two formulas of Theorem~\ref{WelWDVV_thm}, we obtain
\begin{enumerate}[label=($\tau_{1,1}'\arabic*$),leftmargin=*]

\item If $d\!\ge\!2$, $l\!\ge\!1$, and $2d\!-\!l\!\ge\!1$ 
(i.e.~$k\!\ge\!1$), then
\begin{equation*}\begin{split}
{}\hspace{-.8in}
N_{d,l}^{\tau_{1,1}'}=&-(-1)^d2^{l-2}\lr{d}_ld^2N^{\P^1\!\times\P^1}_d 
-\!\!\!\!
\sum_{\begin{subarray}{c}d_0+d'=d\\ d_0,d'\ge1\end{subarray}}\hspace{-.15in}
(-1)^{d'}2^{2d'-2}d_0d'^3\binom{l\!-\!1}{2d'\!-\!1}
N^{\P^1\!\times\P^1}_{d'}N^{\tau_{1,1}'}_{d_0,l-2d'}\\
&+2\hspace{-.2in}
\sum_{\begin{subarray}{c}d_1+d_2=d\\l_1+l_2=l-1\\d_1,d_2\ge1,\,l_1,l_2\ge0\end{subarray}}
\hspace{-.22in}\binom{l\!-\!1}{l_1}\!\!
\Bigg(\!d_1d_2\binom{4d\!-\!2l\!-\!2}{4d_1\!-\!2l_1\!-\!2}\!-\!d_1^2
\binom{4d\!-\!2l\!-\!2}{4d_1\!-\!2l_1\!-\!1}\!\!\!\Bigg)\!
N^{\tau_{1,1}'}_{d_1,l_1}N^{\tau_{1,1}'}_{d_2,l_2}\,.
\end{split}\end{equation*}

\item If $d\!\ge\!2$ and $l\!\ge\!2$, then
\begin{equation*}\begin{split}
{}\hspace{-.8in}
N_{d,l}^{\tau_{1,1}'}=&\!\!\!\!
\sum_{\begin{subarray}{c}d_0+d'=d\\ d_0,d'\ge1\end{subarray}}
\!\!\!(-1)^{d'}2^{2d'-2}d_0d'^2
\Bigg(\!d_0\binom{l\!-\!2}{2d'\!-\!2}-d'\binom{l\!-\!2}{2d'\!-\!1}\!\!\Bigg)\!
N^{\P^1\!\times\P^1}_{d'}N^{\tau_{1,1}'}_{d_0,l-2d'}\\
&+2\hspace{-.2in}
\sum_{\begin{subarray}{c}d_1+d_2=d\\l_1+l_2=l-2\\d_1,d_2\ge1,\,l_1,l_2\ge0\end{subarray}}
\hspace{-.22in}\binom{l\!-\!2}{l_1}\!\Bigg(\!
d_1d_2\binom{4d\!-\!2l\!-\!1}{4d_1\!-\!2l_1\!-\!2}\!-\!
d_1^2\binom{4d\!-\!2l\!-\!1}{4d_1\!-\!2l_1\!-\!1}\!\!\!\Bigg)\!
N^{\tau_{1,1}'}_{d_1,l_1}N^{\tau_{1,1}'}_{d_2,l_2+1}\,.
\end{split}\end{equation*}

\end{enumerate}

\vspace{.10in}

Taking the difference between the above recursions for $N_{d+1,2}^{\tau_{1,1}'}$
with $d\!\ge\!2$, we obtain
$$N_{d,0}^{\tau_{1,1}'}=
\frac{1}{2(d\!-\!1)}
\!\sum_{\begin{subarray}{c}d_1+d_2=d\\ d_1,d_2\ge1\end{subarray}} \!\!
\Bigg(\!d_1\!\binom{4d\!-\!2}{4d_1\!-\!2}\!-\!(d_2\!+\!1)\!\binom{4d\!-\!2}{4d_1\!-\!3}\!\!\Bigg)
N_{d_1,0}^{\tau_{1,1}'}N_{d_2+1,1}^{\tau_{1,1}'}  \qquad\forall\,d\!\ge\!2.$$
The first and third recursions above determine all numbers $N_{d,l}^{\tau_{1,1}'}$ with $d\!\ge\!2$.
The low-degree numbers obtained from these recursions and shown in Table~\ref{WelP1P1a_tbl}
agree with \cite[Cor.~3.18]{Wel07} and \cite[p586]{IKS13}.

\subsection{The product involution}
\label{P1P1a_subs}

Taking $H_1,H_2\!\in\!H^2(\P^1\!\times\!\P^1)$ to be the Poincare duals of~$L_1$ and $L_2$,
respectively, and $B\!=\!aL_1\!+\!bL_2$ 
 in the first formula of Theorem~\ref{WelWDVV_thm},
we obtain
\begin{enumerate}[label=($\tau_{1,1}1\alph*$),leftmargin=*]

\item If $a,b\!\ge\!0$, $l\!\ge\!1$, and $(a\!+\!b)\!-\!l\!\ge\!1$  (i.e.~$k\!\ge\!1$), then
\begin{equation*}\begin{split}
{}\hspace{-1in}
N_{(a,b),l}^{\tau_{1,1}}&=-2^{l-3}ab\,\lr{a,b}_lN^{\P^1\!\times\P^1}_{a/2,b/2}\\
&-\hspace{-.3in}\sum_{\begin{subarray}{c}a_0+2a'=a\\b_0+2b'=b\\
a_0,b_0,a',b'\ge0\\(a_0,b_0),(a',b')\neq(0,0)\end{subarray}} \hspace{-.42in}
2^{2(a'+b')-1}a'b'(a_0b'\!+\!b_0a')\binom{l\!-\!1}{2(a'\!+\!b')\!-\!1}
N^{^{\P^1\!\times\P^1}}_{a',b'}N^{\tau_{1,1}}_{(a_0,b_0),l-2(a'+b')}\\
&+\hspace{-.3in}\sum_{\begin{subarray}{c}a_1+a_2=a\\b_1+b_2=b\\
a_i,b_i\ge0,\,(a_i,b_i)\neq(0,0)\\l_1+l_2=l-1,\,l_1,l_2\ge0\end{subarray}}
\hspace{-.45in}b_1\binom{l\!-\!1}{l_1}
\!\Bigg(\!\!a_2\binom{2(a\!+\!b)\!-\!2l\!-\!2}{2(a_1\!+\!b_1)\!-\!2l_1\!-\!2}\!-\!a_1\binom{2(a\!+\!b)\!-\!2l\!-\!2}{2(a_1\!+\!b_1)\!-\!2l_1\!-\!1}\!\!\!\Bigg)
N^{\tau_{1,1}}_{(a_1,b_1),l_1}N^{\tau_{1,1}}_{(a_2,b_2),l_2}.
\end{split}\end{equation*}

\end{enumerate}
Taking $H_1,H_2\!\in\!H^2(\P^1\!\times\!\P^1)$ to be the Poincare duals of~$L_1$,
$B\!=\!aL_1\!+\!(b\!+\!1)L_2$, and $l\!=\!1$ in
the first formula of Theorem~\ref{WelWDVV_thm}, we obtain
\begin{enumerate}[label=($\tau_{1,1}1\alph*$),leftmargin=*]
\stepcounter{enumi}

\item If $a\!\ge\!2$ and $b\!\ge\!1$, then
\begin{equation*}
{}\hspace{-1in}
N_{(a,b),0}^{\tau_{1,1}}= \frac{1}{2(a\!-\!1)}\hspace{-.1in} 
\sum_{\begin{subarray}{c}a_1+a_2=a\\ b_1+b_2=b+1\\
a_i,b_i\ge1\end{subarray}} \hspace{-.1in}
\Bigg(\!b_1b_2\binom{2(a\!+\!b)\!-\!2}{2(a_1\!+\!b_1)\!-\!2}
\!-\!b_1^2\binom{2(a\!+\!b)\!-\!2}{2(a_1\!+\!b_1)\!-\!1}\!\!\!\Bigg)
N^\phi_{(a_1,b_1),0}N^\phi_{(a_2,b_2),0}\,.
\end{equation*}

\end{enumerate}
These two recursions determine all numbers $N_{(a,b),l}^{\tau_{1,1}}$  with $a\!\ge\!2$.
The low-degree numbers obtained from these recursions and shown in Table~\ref{WelP1P1b_tbl}
agree with \cite[p586]{IKS13}.

\section{Real blowups of $\P^2$}
\label{blowups_sec}

For $k\!\in\!\Z^{\ge0}$, we denote by $\P^2_k$ the blowup of $\P^2$ at 
a generic collection of $k$~points.
For $r,s\!\in\!\Z^{\ge0}$, we denote by~$\tau_{r,s}$ the orientation-reversing involution 
on the blowup
$$\P^2_{r,s}\approx\P^2_{r+2s}$$
of $\P^2$ at $r$ real points and $s$ conjugate pairs of points in~$(\P^2,\tau_2)$
induced by~$\tau_2$.
The fixed locus $(\P^2_{r,s})^{\tau_{r,s}}$ of~$\tau_{r,s}$ 
is the real blowup of $\R\P^2$ at $r$ points. 
Before formulating two recursions for Welschinger's invariants of $(\P^2_{k,s},\tau_{r,s})$
obtained from the first relation of Theorem~\ref{WelWDVV_thm}, 
we recall a recursion of \cite[Thm.~3.6]{GoPa} for the complex GW-invariants of~$\P^2_k$.

\subsection{The complex case}
\label{Cblowups_subs}

Let \hbox{$[k]\!=\!\{1,\ldots,k\}$}.
We denote by $L\!\in\!H_2(\P^2_k;\Z)$ the pull-back of the line class $L$ in~$\P^2$.
For each $i\!\in\![k]$, let $E_i\!\in\!H_2(\P^2_k;\Z)$ be the class of the 
exceptional divisor for the $i$-th blowup point.
The classes~$L$ and $E_1,\ldots,E_k$ freely generate~$H_2(\P^2_k;\Z)$.
If a class  
\BE{Bdcdfn_e}
B_{(d,(c_1,\ldots,c_k))}\equiv dL\!-\!c_1E_1\!-\!\ldots\!-\!c_kE_k \in H_2(\P^2_k;\Z)\EE
is different from all of $E_i$, contains an irreducible curve, and
satisfies $\ell_{\om}(B_{(d,(c_1,\ldots,c_k))})\!\ge\!0$, then
\BE{effcond_e}\begin{split}
d\!\in\!\Z^+, \qquad &c_i\!\in\!\Z^{\ge0}~\forall\,i\!\in\![k], \qquad
\sum_{i\in I}c_i\le d'd~~
\forall\,I\!\subset\!\big\{1,\ldots,k\big\},~2|I|\!\le\!d'(d'\!+\!3),~d'\!\le\!d,\\
&\hspace{.5in} \sum_{i=1}^k c_i\le 3d\!-\!1,  \qquad
\sum_{i=1}^k\!\binom{c_i}{2}\le \binom{d\!-\!1}{2}.
\end{split}\EE
The restrictions on the first line above are obtained by intersecting~\eref{Bdcdfn_e}
with the homology classes 
$$L, \qquad E_i, \qquad\hbox{and}\qquad d'L\!-\!\sum_{i\in I}E_i,$$
respectively; all of these classes contain irreducible curves.
The first condition on the second line is obtained from 
$$c_1(\P^2_k)=3L-\sum_{i=1}^kE_i\,.$$
The last condition in~\eref{effcond_e} follows from arithmetic genus considerations for the image
of a curve in the homology class~\eref{Bdcdfn_e} under the projection to~$\P^2$.\\

For $k\!\in\!\Z^{\ge0}$, let 
$$\cH_k=\Z^+\!\times\!(\Z^{\ge0})^k.$$
For an element $v\!\equiv\!(d,(c_1,\ldots,c_k))$ of $\cH_k$, we define
\begin{gather*}
N^k_v=N^{\P^2_k}_{B_v}, \qquad \ell(v)= 3d\!-\!1-\!\sum_{i=1}^kc_i, \qquad
\wch{v}=d, \\
v0=\big(d,(c_1,\ldots,c_k,0)\big),\quad v1=\big(d,(c_1,\ldots,c_k,1)\!\big),\quad
\cP(v)=\big\{(v_1,v_2)\!\in\cH_k\!\times\!\cH_k\!:v_1\!+\!v_2\!=\!v\big\}.
\end{gather*}
If in addition $i\!\in\![k]$, define 
$$v\!-\!e_i=\big(d,(c_1,\ldots,c_{i-1},c_i\!-\!1,c_{i+1},\ldots,c_k)\!\big),
\qquad \lr{v,e_i}=c_i\,.$$
If $v'\!\equiv\!(d',(c_1',\ldots,c_k'))$ is another element of $\cH_k$, let
$$\lr{v,v'}=B_v\!\cdot_{\P^2_k}\!B_{v'}= dd'-\sum_{i=1}^kc_ic_i'\,.$$

\vspace{.1in}

The (complex) GW-invariants $N^k_v$ vanish unless $v$ satisfies all conditions in~\eref{effcond_e}.
These numbers are preserved by the permutations of the elements 
of the $k$-tuple part of~$v$.
Furthermore, 
$$N_{v0}^{k+1}=N_v^k~\forall\,v\!\in\!\cH_k, ~
N_{v1}^{k+1}=N_v^k~\forall\,v\!\in\!\cH_k~\hbox{s.t.}~\ell(v)\!>\!0, ~
N^{\P^2_k}_{E_i}=1~\forall\,i\!\in\![k],~
N_{(d,())}^0=N_d^{\P^2}~\forall\,d\!\in\!\Z^+.$$
Along with these identities, the next recursion determines all numbers $N^k_v$
with $v\!\in\!\cH_k$ and $k\!\in\!\Z^{\ge0}$. 
\begin{enumerate}[leftmargin=49pt]


\item[{\cite[R(i)]{GoPa}}] If $\ell(v)\!\ge\!0$ and $i\!\in\![k]$, then
\begin{equation*}\begin{split}
{}\hspace{-.8in}
\wch{v}^2\lr{v,e_i}N_v^k&=\big(\wch{v}^2\!-\!(\lr{v,e_i}\!-\!1)^2)N_{v-e_i}^k\\
&\hspace{.2in}+\sum_{(v_1,v_2)\in\cP(v-e_i)}\hspace{-0.3in}
\lr{v_1,v_2}\big(\wch{v}_1\wch{v}_2\lr{v_1,e_i}\lr{v_2,e_i}\!-\!
\wch{v}_1^2\lr{v_2,e_i}^2\big)\binom{\ell(v)}{\ell(v_1)}N_{v_1}^kN_{v_2}^k\,.
\end{split}\end{equation*}

\end{enumerate}

\subsection{The real case}
\label{Rblowups_subs}

For $r,s\!\in\!\Z^{\ge0}$, let 
\BE{Hgendfn_e} E^\R_i\in H_2(\P^2_{r,s};\Z),~~\forall\,i\!\in\![r], 
\qquad\hbox{and}\qquad
E^+_i,E^-_i\in H_2(\P^2_{r,s};\Z),~~\forall\,i\!\in\![s],\EE
be the classes of the exceptional divisors corresponding to the real blowup points
and to the conjugate pairs of blowup points, respectively. 
The subgroup $H_2^{\tau_{r,s}}(\P^2_{r,s})$ of $H_2(\P^2_{r,s};\Z)$ is freely generated 
by the classes
$$L, \qquad E^\R_1,\ldots,E^\R_r, \qquad\hbox{and}\qquad 
E_1^{\C}\!\equiv\!E^+_1\!+\!E^-_1,~\ldots,~E_s^{\C}\!\equiv\!E^+_s\!+\!E^-_s.$$

\vspace{.1in}

For $r,s\in\Z^{\ge0}$, define 
\begin{gather*}
\fd\!:\cH_{r+2s}\lra \cH_{r,s}\!\equiv\!\Z^+\!\times\!(\Z^{\ge0})^r\!\times\!(\Z^{\ge0})^s,\\
\fd\big(d,(c_1,\ldots,c_{r+2s})\big)=
\big(2d,(2c_1,\ldots,2c_r),(c_{r+1}\!+\!c_{r+2},\ldots,c_{r+2s-1}\!+\!c_{r+2s})\big).
\end{gather*}
For an element 
\BE{vdfn_e}v\equiv(d,\a,\b)\equiv \big(d,(a_1,\ldots,a_r),(b_1,\ldots,b_s)\big)\EE
of $\cH_{r,s}$ and $l\!\in\!\Z^{\ge0}$, we define
\begin{gather*}
B_v=dL-\!\sum_{i=1}^ra_iE_i^{\R}-\!\sum_{i=1}^sb_iE_i^{\C}, \quad
N^{r,s}_{v,l}=N^{\tau_{r,s}}_{B_v,l}, \quad 
\ell(v)= 3d\!-\!1-\!\sum_{i=1}^ra_i-2\!\sum_{i=1}^sb_i, \quad
\wch{v}=d, \\
|v|_{\C}=\sum_{i=1}^sb_i,  \qquad
v\!+\!\wch{e}=\big(d\!+\!1,\a,\b\big),\\
\cP_{\R}(v)=\big\{(v_1,v_2)\!\in\cH_{r,s}\!\times\!\cH_{r,s}\!:v_1\!+\!v_2\!=\!v\big\},
\quad
\cP_{\C}(v)=\big\{(v_0,v')\!\in\cH_{r,s}\!\times\!\cH_{r+2s}\!:v_0\!+\!\fd(v')\!=\!v\big\}.
\end{gather*}
If in addition $\ep\!\in\!\{0,1\}$, $i\!\in\![r]$, and $j\!\in\![s]$, let
$$v\ep_{\R}=\big(d,(a_1,\ldots,a_r,\ep),\b\big), \quad
v\ep_{\C}=\big(d,\a,(b_1,\ldots,b_s,\ep)\big), \quad 
\lr{v,e_i^{\R}}=a_i, \quad \lr{v,e_j^{\C}}=2b_j.$$
Define
\begin{gather*}
\lr{\cdot,\cdot}\!:\cH_{r,s}\!\times\!\cH_{r+2s}\lra\Z,\\
\blr{\big(d_0,(a_1,\ldots,a_r),(b_1,\ldots,b_s)\big),\big(d',(c_1,\ldots,c_{r+s})\big)}
=d_0d'\!-\!\sum_{i=1}^r\!a_ic_i\!-\!
\sum_{i=1}^s\!b_i\big(c_{r+2i-1}\!+\!c_{r+2i}\big).
\end{gather*}
For an element $v'\!\equiv\!(d,(c_1,\ldots,c_{r+2s}))$ of $\cH_{r+2s}$, $i\!\in\![r]$,
and $j\!\in\![s]$, let
$$|v'|_{\R}=d\!+\!\sum_{i=1}^rc_i, \quad  \lr{v',e_i^{\R}}=c_i,
\quad  \lr{v',e_j^{\C}}=c_{r+2j-1}\!+\!c_{r+2j}.$$
For $a,b\!\in\!\Z$, let $\de_{a,b}\!=\!1$ if $a\!=\!b$ and $0$ otherwise.\\
 
Welschinger's invariants $N^{r,s}_{v,l}$ with $v\!\in\!\cH_{r,s}$ as above vanish unless
$v$ satisfies all conditions in~\eref{effcond_e} with
$$(c_1,\ldots,c_k)=\big(a_1,\ldots,a_r,b_1,b_1,\ldots,b_s,b_s\big).$$
These numbers are preserved by the permutations of the elements of the tuples~$\a$
and~$\b$.
Furthermore, 
$$N_{v0_{\R},l}^{r+1,s},N_{v0_{\C},l}^{r,s+1}=N_{v,l}^{r,s}
~~\forall\,v\!\in\!\cH_{r,s},\,l\!\in\!\Z^{\ge0},\quad
N^{\tau_{r,s}}_{E_i^{\R},0}=1~~\forall\,i\!\in\![r], \quad
N^{\tau_{r,s}}_{E_i^{\C},0}=0~~\forall\,i\!\in\![s].$$

\vspace{.15in}

Let $\bF\!=\!\R$ and $i\!\in\![r]$ or $\bF\!=\!\C$ and $i\!\in\![s]$.
Taking $H_1,H_2\!\in\!H^2(\P^2_{r+2s})$ to be the Poincare duals of~$L$
and $E^{\bF}_i$, respectively, $B\!=\!B_{v+\wch{e}}$, and $l\!=\!1$ in
the first formula of Theorem~\ref{WelWDVV_thm}, we obtain
\begin{enumerate}[label=($\tau_{r,s}1\alph*$),leftmargin=*]

\item\label{BlDiff2_it} If $\ell(v)\!\ge\!0$ and either $\bF\!=\!\R$ and $i\!\in\![r]$ or
$\bF\!=\!\C$ and $i\!\in\![s]$, then 
\begin{equation*}\begin{split}
{}\hspace{-.2in}
\big(\wch{v}\!-\!\ell(v)\big)\lr{v,e^\bF_i}
N^{r,s}_{v,0}
=&-\frac{(-1)^{|v|_{\C}}\de_{\ell(v),0}}{4}\lr{v,e^\bF_i}\big(\wch{v}\!+\!1\big)
\hspace{-.18in}
\sum_{\begin{subarray}{c}v'\in\cH_{r+2s}\\ \fd(v')=v+\wch{e}\end{subarray}}
\hspace{-.2in}N^{r+2s}_{v'}\\
&-\hspace{-0.2in}\sum_{\begin{subarray}{c}(v_0,v')\in\cP_\C(v+\wch{e})\\ \ell(v')=0\end{subarray}}
\hspace{-0.3in}(-1)^{|v'|_{\R}}\lr{v_0,v'}\wch{v}'\lr{v',e^\bF_i}
N^{r+2s}_{v'}N^{r,s}_{v_0,0}\,\\
&+\hspace{-0.35in}\sum_{\begin{subarray}{c}(v_1,v_2)\in\cP_\R(v+\wch{e})\\
v_1,v_2\neq v\end{subarray}}\hspace{-0.2in}\wch{v}_1\Bigg(\!\!
\lr{v_2,e^\bF_i}\binom{\ell(v)}{\ell(v_1)\!-\!1}
\!-\!\lr{v_1,e^\bF_i}\binom{\ell(v)}{\ell(v_1)}\!\!\!\Bigg)
\!N^{r,s}_{v_1,0}N^{r,s}_{v_2,0}.
\end{split}\end{equation*}
\end{enumerate}
Taking $H_1,H_2\!\in\!H^2(\P^2_{r+2s})$ to be the Poincare duals of $L\!-\!E^{\R}_i$
with $i\!\in\![r]$ instead, $B\!=\!B_{v-2e_i^{\R}}$, and $l\!=\!1$ in
the first formula of Theorem~\ref{WelWDVV_thm}, we obtain
\begin{enumerate}[label=($\tau_{r,s}1\alph*$),leftmargin=*]

\setcounter{enumi}{1}

\item\label{BlDiff2b_it}
If $\ell(v)\!\ge\!1$, $i\!\in\![r]$, and $\lr{v,e^{\R}_i}\!\ge\!2$, then 
\begin{equation*}\begin{split}
&{}\hspace{-.7in}
\lr{v,e^{\R}_i}N^{r,s}_{v,0}
=N^{r,s}_{v-e^{\R}_i,0}
+\frac{(-1)^{|v|_{\C}}\de_{\ell(v),1}}{4}\big(\wch{v}\!-\!\lr{v,e^{\R}_i}\!+\!2\big)^2
\hspace{-.18in}
\sum_{\begin{subarray}{c}v'\in\cH_{r+2s}\\ \fd(v')=v-2e_i^{\R}\end{subarray}}
\hspace{-.25in}N^{r+2s}_{v'}\\
&\hspace{-0.4in}
+\hspace{-0.2in}\sum_{\begin{subarray}{c}(v_0,v')\in\cP_\C(v-2e_i^{\R})\\ \ell(v')=0\end{subarray}}
\hspace{-0.37in}(-1)^{|v'|_{\R}}\lr{v_0,v'}\big(\wch{v}'\!-\!\lr{v',e^\R_i}\big)^{\!2}
N^{r+2s}_{v'}N^{r,s}_{v_0,0}\,\\
&\hspace{-0.4in}
-\hspace{-0.2in}\sum_{(v_1,v_2)\in\cP_\R(v-2e_i^{\R})}\hspace{-0.42in}
\big(\wch{v}_1\!-\!\lr{v_1,e^{\R}_i}\!\big)\Bigg(\!\!
\big(\wch{v}_2\!-\!\lr{v_2,e^{\R}_i}\big)\!\binom{\ell(v)\!-\!1}{\ell(v_1)\!-\!1}
\!-\!\big(\wch{v}_1\!-\!\lr{v_1,e^{\R}_i}\big)\!\binom{\ell(v)\!-\!1}{\ell(v_1)}\!\!\!\Bigg)
\!N^{r,s}_{v_1,0}N^{r,s}_{v_2,0}.
\end{split}\end{equation*}
\end{enumerate}

\vspace{.15in}

Let $v\!\in\!\cH_{r,s}$ be as in~\eref{vdfn_e} with $\ell(v)\!\ge\!0$.
If $a_i\!=\!0$ for all $i\!\in\![r]$ and $b_j\!=\!0$ for all $j\!\in\![s]$,
then the numbers
$$N_{v,l}^{r,s}=N_{(d,(),()),l}^{0,0}=N_{d,l}^{\tau_2}$$
are computable from the recursions for Welschinger's invariants of $(\P^2,\tau_2)$ 
provided by Section~\ref{P2_sec}.
Since
$$N_{v1_{\C},l}^{r,s+1}=N_{v,l+1}^{r,s} \qquad\forall\,v\!\in\!\cH_{r,s},\,l\!\in\!\Z^{\ge0}\,,$$
it suffices in general to determine Welschinger's invariants $N_{v,0}^{r,s}$ with
real point insertions only.
Since 
$$N_{v'1_{\R},l}^{r,s}=N_{v',l}^{r-1,s} \qquad\forall\,v'\!\in\!\cH_{r-1,s},\,l\!\in\!\Z^{\ge0}
~~\hbox{with}~~\ell(v')\!>\!2l,\,$$
we can also assume that $a_i\!\ge\!2$ for all $i\!\in\![r]$.
If $a_i\!\ge\!2$ for some $i\!\in\![r]$, 
$N_{v,0}^{r,s}$ is reduced to the numbers $N_{v',0}^{r,s}$ with $|v'|\!<\!|v|$
by the relation~\ref{BlDiff2b_it} if $\ell(v)\!>\!0$ and  
by the $\bF\!=\!\R$ case of~\ref{BlDiff2_it} if $\ell(v)\!=\!0$.
If $a_i\!=\!0$ for all $i\!\in\![r]$ and $b_j\!\neq\!0$ for some $j\!\in\![s]$,
then the $(\bF,i)\!=\!(\C,j)$ case of~\ref{BlDiff2_it}
expresses $N_{v,0}^{r,s}$ in terms of the numbers $N_{v',0}^{r,s}$ with $|v'|\!<\!|v|$.\\

Welschinger's invariants $N_{v,l}^{r,s}$ of $(\P^2_{r,s},\tau_{r,s})$
obtained from the above recursions and 
shown in Tables~\ref{WelP2bl1R_tbl}-\ref{WelP2bl3C_tbl}
agree with \cite[Table~4]{B15}, \cite[Corollary~3.2]{BP13},
\cite[pp585/6]{IKS13}, and \cite[Example~17]{IKS13a}.
Other numbers obtained from these recursions agree with \cite[Tables~6,9]{B15}
and \cite[Section~2.3]{IKS15}.

\section{WDVV-type relations for real symplectic sixfolds}
\label{genWel3_sec}

We now take $(X,\om,\phi)$ to be a compact real symplectic sixfold
with connected orientable fixed locus~$X^{\phi}$.
An \sf{automorphism} of $(X,\om,\phi)$ is a diffeomorphism~$\psi$ of~$X$ 
such~that 
$$\psi^*\om=\om \qquad\hbox{and}\qquad \psi\!\circ\!\phi=\phi\!\circ\!\psi.$$
We call such an automorphism an \sf{averager} if $\psi$ is an involution
which acts trivially on $H_2^{\phi}(X)$ and restricts to 
an orientation-reversing diffeomorphism of~$X^{\phi}$. 
As explained in Section~2.5 of~\cite{RealWDVV3}, 
these conditions imply that the natural homomorphisms
\BE{aveG_e0}\io_*\!:H_2\big(X\!-\!X^{\phi};\Z\big)\lra H_2(X;\Z) 
\quad\hbox{and}\quad
r\!:H^4\big(X,X^{\phi}\big)\lra H^4(X)\EE
restrict to isomorphisms
\begin{equation*}\begin{split}
\io_*^{\psi}\!:\big\{B\!\in\!H_2(X\!-\!\wch{X}^{\phi};\Z)\!:
\phi_*B\!=\!-B,\,\psi_*B\!=\!B\big\} &\stackrel{\approx}{\lra} H_2^{\phi}(X)
\qquad\hbox{and}\\
r_{\psi}\!:\big\{\mu\!\in\!H^4(X,\wch{X}^{\phi})\!:
\phi^*\mu\!=\!\mu,\,\psi^*\mu\!=\!\mu\big\}&\stackrel{\approx}{\lra} H^4(X)^{\phi}_+\,,
\end{split}\end{equation*}
respectively.
For a homogeneous element $\mu$ of $H^{2*}(X)$, let
\BE{wtmudfn_e}\wt\mu=\begin{cases}r_{\psi}^{-1}(\mu), &\hbox{if}~\mu\!\in\!H^4(X)^{\phi}_+;\\
0, &\hbox{if}~\mu\!\in\!H^4(X)^{\phi}_-;\\
\mu,&\hbox{if}~\mu\!\not\in\!H^4(X).
\end{cases}\EE

\vspace{.15in}

Since $X^{\phi}$ is an orientable three-dimensional manifold,
its tangent bundle is trivializable and thus admits a Spin-structure~$\fs$ for any choice of
orientation~$\fo$ in~$X^{\phi}$.
We call such a pair $\os\!\equiv\!(\fo,\fs)$ an \sf{OSpin-structure on~$X^{\phi}$}.
There is an associated OSpin-structure $\ov\os\!\equiv\!(\ov\fo,\ov\fs)$ 
for the opposite orientation~$\ov\fo$ of~$\fo$ on~$X^{\phi}$;
see the SpinPin~3 property in Section~1.2 of~\cite{SpinPin}.
If $\os'$ is another OSpin-structure on~$X^{\phi}$ and $b\!\in\!H_1(X^{\phi};\Z_2)$,
we write $\os|_b\!=\!\os'|_b$ if $\al^*\os\!=\!\al^*\os'$ for any loop $\al\!:S^1\!\lra\!X^{\phi}$
representing~$b$.
The group $H^1(X^{\phi};\Z_2)$ acts on the set of OSpin-structures
on~$X^{\phi}$ associated with a fixed orientation~$\fo$ freely and transitively;
see the SpinPin~2 property in Section~1.2 in~\cite{SpinPin}.
If $\os'\!=\!\eta\!\cdot\!\os$, then $\os|_b\!=\!\os'|_b$ if and only if $\lr{\eta,b}\!=\!0$.
For $B\!\in\!H_2(X;\Z)$ and a tuple $(\mu_1,\ldots,\mu_l)$ of homogeneous elements 
of $H^{2*}(X)$ and $H^{2*}(X,X^{\phi})$, let 
\BE{dimcond3_e} 
k= k_B(\mu_1,\ldots,\mu_l)
\equiv\frac12\Big(\ell_{\om}(B)\!+\!2l-\sum_{i=1}^l\deg\mu_i\Big).\EE
Under certain conditions on $B$ and $\mu_1,\ldots,\mu_l$,
an OSpin-structure $\os$ on~$X^{\phi}$ determines an \sf{open GW-invariant} 
\BE{RGWdfn_e}\blr{\mu_1,\ldots,\mu_l}_B^{\phi,\os}\in\Q\EE
of $(X,\om,\phi)$
enumerating real irreducible degree~$B$ $J$-holomorphic curves $C\!\subset\!X$ 
that meet~$X^{\phi}$ and
pass through generic representatives for the Poincare duals of $\mu_1,\ldots,\mu_l$
and through~$k$ points in~$X^{\phi}$.
These conditions are recalled in the next paragraph.
If such curves exist, then $\ell_{\om}(B)$ is even and thus $k\!\in\!\Z$.
The number~\eref{RGWdfn_e} is defined to be~0 if $k\!<\!0$.\\

Invariant signed counts~\eref{RGWdfn_e} were first defined in~\cite{Wel6} under the assumptions 
that 
\BE{taucond_e1}\ell_{\om}(B)>0, \qquad \mu_i\!\in\!H^2(X)\!\cup\!H^6(X)~~\forall\,i,\EE
i.e.~each $\mu_i$ represents a Poincare dual of a ``complex" hypersurface or a point, 
and either 
\BE{noS2bub_e} k\!>\!0 \qquad\hbox{or}\qquad
B\not\in \fd\big(H_2(X;\Z)\!\big) \subset H_2^{\phi}(X).\EE
The interpretation of these counts in terms of $J$-holomorphic maps from disks 
in~\cite{Jake} dropped the first restriction in~\eref{taucond_e1} and 
allowed insertions $\mu_i\!\in\!H^4(X,X^{\phi})$.\\

As shown in~\cite{RealWDVV3}, the restriction~\eref{noS2bub_e} can be dropped 
if $(X,\om,\phi)$ admits an averager~$\psi$.
With the notation as in~\eref{wtmudfn_e}, we then define 
\BE{Gnumsdfn_e}\blr{\mu_1,\ldots,\mu_l}_{B;\psi}^{\phi,\os}
=\blr{\wt\mu_1,\ldots,\wt\mu_l}_B^{\phi,\os}\EE
for all $\mu_1,\ldots,\mu_l\!\in\!H^{2*}(X)$.
These numbers satisfy the usual divisor relation, i.e.
\BE{Rdivrel_e}\blr{\mu,\mu_1,\ldots,\mu_l}_{B;\psi}^{\phi,\os}
=\lr{\mu,B}\blr{\mu_1,\ldots,\mu_l}_{B;\psi}^{\phi,\os}
\quad\forall\,\mu\!\in\!H^2(X)\,.\EE
In general, the numbers~\eref{Gnumsdfn_e}
 depend on the $\psi$-invariant subspace of $H^4(X,X^{\phi})$
or equivalently on the $\psi$-invariant subspace of $H_2(X\!-\!X^{\phi})$.
However, they do not depend on the choice of~$\psi$ 
which acts trivially on~$H^2(X)$ if the subspace of~$H^4(X)$ spanned by 
the cup products of the elements of~$H^2(X)$ contains~$H^4(X)_+^{\phi}$.
This is the case in all three examples of Sections~\ref{P3_sec} and~\ref{P1P1P1_sec}.\\

We denote by 
$$\prt_{X^{\phi};\Z_2}\!: H_2\big(X,X^{\phi};\Z\big)\lra H_1\big(X^{\phi};\Z\big)
\lra H_1\big(X^{\phi};\Z_2\big)$$
the composition of the boundary homomorphism of the relative exact sequence for the pair 
$(X,X^{\phi})$ with the mod~2 reduction of the coefficients.
Let 
$$\fd_{X^{\phi}}\!:H_2\big(X,X^{\phi};\Z\big)\lra H_2(X;\Z)$$
be the homomorphism obtained by gluing each continuous map 
$f\!:(\Si,\prt \Si)\!\lra\!(X,X^{\phi})$ from a bordered surface with
the map~$\phi\!\circ\!f$ from~$\Si$ with the opposite orientation along~$\prt\Si$.
By Proposition~1.2 and Remark~1.3 in~\cite{RealWDVV3}, 
the numbers~\eref{Gnumsdfn_e} satisfy the following vanishing property.

\begin{thm}\label{WelWDVV3a_thm}
Suppose $(X,\om,\phi)$ is a compact real symplectic sixfold with connected fixed locus~$X^{\phi}$,
$\psi$ is an averager for $(X,\om,\phi)$,
$\os$ is an OSpin-structure on~$X^{\phi}$, 
$B\!\in\!H_2^{\phi}(X)$, $\mu_1,\ldots,\mu_l\!\in\!H^{2*}(X)$,
and $k\!\in\!\Z$ is as in~\eref{dimcond3_e}.
The numbers~\eref{Gnumsdfn_e} vanish if any of the following conditions holds:
\begin{enumerate}[label=(\alph*),leftmargin=*]

\item $\mu_i\!\in\!H^2(X)^{\phi}_+\!\oplus\!H^4(X)^{\phi}_-$ for some~$i$;

\item $k\!\in\!2\Z$ and $\prt_{X^{\phi};\Z_2}(B')\!=\!0$ for every $B'\!\in\!\fd_{X^{\phi}}^{-1}(B)$;

\item $k\!\in\!2\Z$ and $(\psi^*\os)|_b\!=\!\ov\os|_b$ for every
$b\!\in\!\prt_{X^{\phi};\Z_2}(\fd_{X^{\phi}}^{-1}(B)\!)$;

\item $k\!\not\in\!2\Z$ and $(\psi^*\os)|_b\!\neq\!\ov\os|_b$ for every
$b\!\in\!\prt_{X^{\phi};\Z_2}(\fd_{X^{\phi}}^{-1}(B)\!)$.

\end{enumerate}
\end{thm}

\vspace{.1in}

Choose a basis $\mu_1^{\st},\ldots,\mu_N^{\st}$ for 
$$H^0(X)\!\oplus\!H^2(X)^{\phi}_-\!\oplus\!H^4(X)^{\phi}_+\!\oplus\!H^6(X)$$
consisting of homogeneous elements.
Let $(g_{ij})_{i,j}$ be the $N\!\times\!N$-matrix given~by
$$g_{ij}=\blr{\mu_i^{\st}\mu_j^{\st},[X]}$$
and $(g^{ij})_{ij}$ be its inverse.
For $\mu_1,\ldots,\mu_l\!\in\!H^{2*}(X)$, define
$$\blr{\mu_1,\ldots,\mu_l}_{0;\psi}^{\phi,\os}
=\begin{cases}\lr{\mu_1,[\pt]},&\hbox{if}\,l\!=\!1;\\
0,&\hbox{otherwise}.\end{cases}$$
We denote by $\lr{\mu_1,\ldots,\mu_l}_B^X\!\in\!\Q$
the (complex) GW-invariant of $(X,\om)$ enumerating rational degree~$B$ 
$J$-holomorphic curves $C\!\subset\!X$ through generic representatives of 
the Poincare duals of $\mu_1,\ldots,\mu_l$.
If in addition $I\!\subset\!\{1,2,\ldots,l\}$, 
let $\mu_I$ be the $|I|$-tuple consisting of the entries of $\mu\!\equiv\!(\mu_1,\ldots,\mu_l)$
indexed by~$I$.
For $i,j\!=\!1,2,\ldots,l$, we define
\begin{alignat*}{2}
\cP(l)&=\big\{(I,J)\!:\{1,2,\ldots,l\}\!=\!I\!\sqcup\!J,~1\!\in\!I\big\},&\quad
\cP_{i;}(l)&=\big\{(I,J)\!\in\!\cP(l)\!:i\!\in\!I\big\},\\
\cP_{;j}(l)&=\big\{(I,J)\!\in\!\cP(l):j\!\in\!J\big\},&\quad
\cP_{i;j}(l)&=\cP_{i;}(l)\!\cap\!\cP_{;j}(l).
\end{alignat*}
Let $[N]\!=\!\{1,2,\ldots,N\}$.

\begin{thm}\label{WelWDVV3_thm}
Let $(X,\om,\phi)$, $\psi$, $\os$, $B$, and 
$\mu\!\equiv\!(\mu_1,\ldots,\mu_l)$ be as in Theorem~\ref{WelWDVV3a_thm}
with 
$$k\equiv\frac12\Big(\ell_{\om}(B)\!+\!2l-\sum_{i=1}^l\deg\mu_i\Big)\!-\!1\ge0.$$
\begin{enumerate}[label=($\R$WDVV\arabic*),leftmargin=*]

\item\label{Wel12rec3_it} If $l\!\ge\!2$ and $k\!\ge\!1$, then
\begin{equation*}\begin{split}
{}\hspace{-1in}
\sum_{\begin{subarray}{c}B_0,B'\in H_2(X)\\ B_0+\fd(B')=B\\ 
(I,J)\in\cP_{2;}(l)\end{subarray}}
\!\!\!\!\!\!\!2^{|I|-2}\!\!\!\!
\sum_{i,j\in[N]}\!\!\!\!
\blr{\mu_I,\mu^\st_i}^X_{B'}g^{ij}\blr{\mu^\st_j,\mu_J}^{\phi,\os}_{B_0;\psi}
&+\sum_{\begin{subarray}{c}B_1,B_2\in H_2^{\phi}(X)\\ B_1+B_2=B\\  
(I,J)\in\cP_{2;}(l)\end{subarray}}
\hspace{-.1in}\binom{k\!-\!1}{k_{B_1}(\mu_I)}
\blr{\mu_I}^{\phi,\os}_{B_1;\psi}\blr{\mu_J}^{\phi,\os}_{B_2;\psi}\\
&=\sum_{\begin{subarray}{c}B_1,B_2\in H_2^{\phi}(X)\\ B_1+B_2=B\\ 
(I,J)\in\cP_{;2}(l)\end{subarray}}
\hspace{-.1in}
\binom{k\!-\!1}{k_{B_1}(\mu_I)\!-\!1}
\blr{\mu_I}^{\phi,\os}_{B_1;\psi}
\blr{\mu_J}^{\phi,\os}_{B_2;\psi}\,.
\end{split}\end{equation*}

\item\label{Wel03rec3_it} If $l\!\ge\!3$, then
\begin{equation*}\begin{split}
&\hspace{-1in}
\sum_{\begin{subarray}{c}B_0,B'\in H_2(X)\\ B_0+\fd(B')=B\\ 
(I,J)\in\cP_{2;3}(l)\end{subarray}}
\!\!\!\!\!\!\!2^{|I|-2}\!\!\!\!
\sum_{i,j\in[N]}\!\!\!\!
\blr{\mu_I,\mu^\st_i}^X_{B'}g^{ij}
\blr{\mu^{\st}_j,\mu_J}^{\phi,\os}_{B_0;\psi}
+\sum_{\begin{subarray}{c}B_1,B_2\in H_2^{\phi}(X)\\ B_1+B_2=B\\ 
(I,J)\in\cP_{2;3}(l)\end{subarray}}
\hspace{-.1in}
\binom{k}{k_{B_1}(\mu_I)}
\blr{\mu_I}^{\phi,\os}_{B_1;\psi}\blr{\mu_J}^{\phi,\os}_{B_2;\psi}\\
&\hspace{-.8in}
=\sum_{\begin{subarray}{c}B_0,B'\in H_2(X)\\ B_0+\fd(B')=B\\ 
(I,J)\in\cP_{3;2}(l)\end{subarray}}
\!\!\!\!\!\!\!2^{|I|-2}\!\!\!\!
\sum_{i,j\in[N]}\!\!\!\!
\blr{\mu_I,\mu^\st_i}^X_{B'}g^{ij}
\blr{\mu^{\st}_j,\mu_J}^{\phi,\os}_{B_0;\psi}
+\sum_{\begin{subarray}{c}B_1,B_2\in H_2^{\phi}(X)\\ B_1+B_2=B\\ 
(I,J)\in\cP_{3;2}(l)\end{subarray}}
\hspace{-.1in}
\binom{k}{k_{B_1}(\mu_I)}
\blr{\mu_I}^{\phi,\os}_{B_1;\psi}\blr{\mu_J}^{\phi,\os}_{B_2;\psi}\,.
\end{split}\end{equation*}

\end{enumerate}
\end{thm}

\vspace{.2in}

Let $\cN X^{\phi}$ denote the normal bundle of $X^{\phi}$ in~$X$.
An orientation~$\fo$ on~$X^{\phi}$ determines an orientation 
on a unit sphere $S(\cN_pX^{\phi})$ in the fiber of~$\cN X^{\phi}$ 
over any $p\!\in\!X^{\phi}$.
As explained in Section~2.5 of~\cite{RealWDVV3},
the first homomorphism in~\eref{aveG_e0} is surjective with
the kernel generated by the homology class $[S(\cN_pX^{\phi})]_{X-X^{\phi}}$
of~$S(\cN_pX^{\phi})$ in $X\!-\!X^{\phi}$;
the second homomorphism in~\eref{aveG_e0} is related to the first by the Poincare Duality.
According to Proposition~2.1 in~\cite{RealWDVV3},
\BE{SNXphi_e}\blr{\mu_1,\ldots,\mu_l,
\PD_{X,X^{\phi}}\big([S(\cN_pX^{\phi})]_{X-X^{\phi}}\big)}_B^{\phi,\os}
=2\blr{\mu_1,\ldots,\mu_l}_B^{\phi,\os}\EE
for all $B\!\in\!H_2^{\phi}(X)$ and 
$\mu_i\!\in\!H^2(X)\!\oplus\!H^4(X,X^{\phi})\!\oplus\!H^6(X)$,
i.e.~an insertion of $S(\cN_pX^{\phi})$ can be traded for a real point insertion.
Thus, the invariants~\eref{Gnumsdfn_e} with $\mu_i\!\in\!H^{2*}(X)$ determine 
the invariants~\eref{RGWdfn_e} with $\mu_i\!\in\!H^2(X)\!\oplus\!H^4(X,X^{\phi})\!\oplus\!H^6(X)$.

\begin{rmk}\label{WelWDVV3_rmk}
Theorems~\ref{WelWDVV3a_thm} and~\ref{WelWDVV3_thm} extend to 
compact real symplectic sixfolds $(X,\om,\phi)$ with disconnected fixed loci~$X^{\phi}$
admitting finite groups of symmetries with certain properties
and apply with finer notions of the curve degree~$B$;
see Proposition~1.2, Remark~1.3, and Theorem~1.4 in~\cite{RealWDVV3}.
This makes no difference for the examples of Sections~\ref{P3_sec} and~\ref{P1P1P1_sec}.
\end{rmk}

\section{The projective space $\P^3$}
\label{P3_sec}

The fixed locus of the standard conjugation 
$$\tau_3\!:\P^3\lra\P^3, \qquad \tau_3\big([Z_0,Z_1,Z_2,Z_3]\big)
=\big[\ov{Z_0},\ov{Z_1},\ov{Z_2},\ov{Z_3}\big],$$
is the real projective space~$\R\P^3$.
The group \hbox{$H_2^{\tau_3}(\P^3)\!=\!H_2(\P^3;\Z)$}
is identified with~$\Z$ via the standard generator $L\!=\![\P^1]$.
An averager~$\psi$ in this case is a reflection about a $\tau_3$-invariant 
complex hyperplane such~as
$$\psi_3\!:\P^3\lra\P^3, \qquad  
\psi_3\big([Z_0,Z_1,Z_2,Z_3]\big)=\big[Z_0,Z_1,Z_2,-Z_3\big].$$
Each of the two orientations on~$\R\P^3$ is compatible with two Spin-structures.
A canonical OSpin-structure~$\os_0$ on~$\R\P^3$ is specified in Section~2.2 
of~\cite{RGWsIII}.   
It is straightforward to see~that
\BE{P3_e3}\begin{aligned}
\prt_{\R\P^3;\Z_2}(B')&=0 &\quad&\hbox{if}~\fd_{\R\P^3}(B')\in2\Z L 
\qquad\hbox{and}\\
\psi^*\os_0\big|_{\prt_{\R\P^3;\Z_2}(B')}&\neq
\ov{\os_0}\big|_{\prt_{\R\P^3;\Z_2}(B')} &\quad&\hbox{if}~
\fd_{\R\P^3}(B')\not\in2\Z L.
\end{aligned}\EE
By the naturality of the $H^1(\R\P^3;\Z_2)$-action on the set of OSpin-structures
on~$\R\P^3$, the last inequality holds for all OSpin-structures on~$\R\P^3$.\\

Let $H\!\in\!H^2(\P^3)$ denote the hyperplane class.
For $d\!\in\!\Z^+$ and a tuple $\m\!\equiv\!(m_1,\ldots,m_l)$ of nonnegative integers,
define
\begin{alignat*}{2}
k_d(\m)&= 2d\!+\!l\!-\!\sum_{i=1}^lm_i,  &\quad
\lr{d}_{\m}&=\begin{cases}1,&\hbox{if}~d\!\in\!2\Z,\,k_d(\m)\!=\!2;\\
0,&\hbox{otherwise};
\end{cases}\\
\blr{\m}_{\!d}&=\blr{H^{m_1},\ldots,H^{m_l}}_{\!dL}^{\!\P^3}, &\quad
\blr{\m}_{\!d}^{\!\tau_3}&=\blr{H^{m_1},\ldots,H^{m_l}}_{\!dL;\psi_3}^{\!\tau_3,\os_0}\,.
\end{alignat*}
If in addition $I\!\subset\![l]$, 
let $\m_I$ be the $|I|$-tuple consisting of the entries of~$\m$ indexed by~$I$.\\

Since there is a unique line~$\P^1$ passing through every pair of points in~$\P^3$
and the orientation conventions of~\cite{RealWDVV3} yield the opposite sign
of Example~6.3 in~\cite{Teh}, 
$$ \blr{3,3}_{\!1}=1, \qquad \blr{3}_{\!1}^{\!\tau_3}=-1\,.$$
By~\eref{P3_e3} and Theorem~\ref{WelWDVV3a_thm},
$$\blr{\m}_{\!d}^{\!\tau_3}=0 \qquad\hbox{if}\quad d\equiv k_d(\m)~~\hbox{mod}~2.$$
The $\P^3$ case of Theorem~10.4 in~\cite{RT} can be written~as
\begin{equation*}\begin{split}
\blr{m_1,m_2,m_3\!+\!1,m_4,\ldots,m_l}_{\!d}&=\blr{m_1,m_2\!+\!1,m_3,m_4,\ldots,m_l}_{\!d}
\!+\!d\blr{m_1\!+\!m_3,m_2,m_4,\ldots,m_l}_{\!d}\\
&\qquad+
\sum_{\begin{subarray}{c}d_1,d_2\in\Z^+\\ d_1+d_2=d\end{subarray}}
\!\!\!\Bigg\{\!\!\sum_{(I,J)\in\cP_{3;2}(l)}\!-\!
\sum_{(I,J)\in\cP_{2;3}(l)}\!\!\Bigg\}
d_2\!\!\!\!
\sum_{\begin{subarray}{c}i,j\in\Z^+\\ i+j=3\end{subarray}}\!\!\!\!
\blr{\m_I,i}_{\!d_1}\blr{j,\m_J}_{\!d_2}
\end{split}\end{equation*}
whenever $m_1,m_2\!\ge\!2$.
Along with the divisor relation, the above identity with 
$$m_1\ge m_2\ge m_3\!+\!1 \ge m_4\ge \ldots\ge m_l\ge2$$
recursively determines all complex GW-invariants $\lr{\m}_{\!d}$ of~$\P^3$.\\

For the purposes of applying Theorem~\ref{WelWDVV3_thm} in this case,
we take the basis $\{\mu_i^{\st}\}$ for~$H^*(\P^3)$ to be~$1,H,H^2,H^3$.
Taking $\mu_2\!=\!H$ in~\ref{Wel12rec3_it} and $\mu_3\!=\!H$ in~\ref{Wel03rec3_it},
we obtain
\begin{enumerate}[label=($\tau_3\arabic*$),leftmargin=*]

\item If $d\!\ge\!1$, $l\!\ge\!1$, and $\m\!\equiv\!(m_1,\ldots,m_l)\!\in\!(\Z^+)^l$
with $k_d(\m)\!\ge\!2$, then
\begin{equation*}\begin{split}
\hspace{-.4in}\blr{m_1\!+\!1,m_2,\ldots,m_l}_{\!d}^{\!\tau_3}
&=-2^{l-2}\lr{d}_{\m}d\blr{m_1,m_2,\ldots,m_l,3}_{\!d/2}\\
&\quad-\!\!\!
\sum_{\begin{subarray}{c}d_0,d'\in\Z^+\\ d_0+2d'=d\end{subarray}}
d'\!\!\!\!\!\!\!\sum_{(I,J)\in\cP(l)}
\!\!\!\!\!\!2^{|I|-1}\!\!\!\!\!\sum_{\begin{subarray}{c}i,j\in\Z^+\\ i+j=3\end{subarray}}\!\!\!\!
\blr{\m_I,i}_{\!d'}\blr{j,\m_J}_{\!d_0}^{\!\tau_3}\\
&\quad+\!\!\!\sum_{\begin{subarray}{c}d_1,d_2\in\Z^+\\ d_1+d_2=d\end{subarray}}
\sum_{(I,J)\in\cP(l)}\!\!\!
\Bigg(\!d_2\!\binom{k_d(\m)\!-\!2}{k_{d_1}(\m_I)\!-\!1}\!-\!
d_1\!\binom{k_d(\m)\!-\!2}{k_{d_1}(\m_I)}\!\!\Bigg)
\blr{\m_I}_{\!d_1}^{\!\tau_3}\blr{\m_J}_{\!d_2}^{\!\tau_3}\,.
\end{split}\end{equation*}

\item If $d\!\ge\!1$, $l\!\ge\!2$, and $\m\!\equiv\!(m_1,\ldots,m_l)\!\in\!(\Z^+)^l$
with $k_d(\m)\!\ge\!1$, then
\begin{equation*}\begin{split}
&{}\hspace{-.4in}\blr{m_1,m_2\!+\!1,m_3,\ldots,m_l}_{\!d}^{\!\tau_3}
=\blr{m_1\!+\!1,m_2,m_3,\ldots,m_l}_{\!d}^{\!\tau_3}\\
&+\!\!\!
\sum_{\begin{subarray}{c}d_0,d'\in\Z^+\\ d_0+2d'=d\end{subarray}}
d'\!\!\!\!\!\!\!\sum_{(I,J)\in\cP_{;2}(l)}
\sum_{\begin{subarray}{c}i,j\in\Z^+\\ i+j=3\end{subarray}}\!\!\!
\Big(2^{|I|-1}\blr{\m_I,i}_{\!d'}\blr{j,\m_J}_{\!d_0}^{\!\tau_3}
-2^{|J|-1}\blr{\m_J,j}_{\!d'}\blr{i,\m_I}_{\!d_0}^{\!\tau_3}\!\Big)\\
&+\!\!\!\sum_{\begin{subarray}{c}d_1,d_2\in\Z^+\\ d_1+d_2=d\end{subarray}}
\sum_{(I,J)\in\cP_{;2}(l)}\!\!\!
\Bigg(\!\!d_1\binom{k_d(\m)\!-\!1}{k_{d_1}(\m_I)}\!-\!
d_2\binom{k_d(\m)\!-\!1}{k_{d_2}(\m_J)}\!\!\Bigg)
\blr{\m_I}_{\!d_1}^{\!\tau_3}\blr{\m_J}_{\!d_2}^{\!\tau_3}\,.
\end{split}\end{equation*}

\end{enumerate}

\vspace{.15in}

From the two equations above and the real divisor relation~\eref{Rdivrel_e}, we get
\begin{gather*}
\lr{2}_{\!1}^{\!\tau_3}=0, \qquad \lr{2,2}_{\!1}^{\!\tau_3}=\lr{3,1}_{\!1}^{\!\tau_3}=-1,\\
\lr{3,2}_{\!2}^{\!\tau_3}=-2\lr{2,2,3}_{\!1}\!-\!\lr{2,2}_{\!1}^{\!\tau_3}\lr{}^{\!\tau_3}_{\!1}
\!+\!\lr{2}^{\!\tau_3}_{\!1}\lr{2}^{\!\tau_3}_{\!1}, \quad
\lr{3,2}_{\!2}^{\!\tau_3}=\lr{3}_{\!1}^{\!\tau_3}\lr{1}^{\!\tau_3}_{\!1}, \quad
\lr{}^{\!\tau_3}_{\!1}=1.
\end{gather*}
Taking the difference between the two equations for $\lr{3,2}_{d+1}^{\tau_3}$ 
with $d\!\ge\!2$, we obtain
$$\blr{}_{\!d}^{\!\tau_3}=\frac1{d\!+\!1}
\!\!\!\sum_{\begin{subarray}{c}d_1,d_2\in\Z^+\\ d_1+d_2=d\end{subarray}}
\!\!\!
\Bigg(\!d_2\!\binom{2d\!-\!1}{2d_1}\lr{3}_{\!d_1+1}^{\!\tau_3}\lr{}_{\!d_2}^{\!\tau_3}
\!+\!\frac{2d_2\!-\!d_1\!-\!1}{2d_2\!-\!1}\!
\binom{2d\!-\!2}{2d_1}\!
\Big(\lr{2,2}_{\!d_1+1}^{\!\tau_3}\lr{}_{\!d_2}^{\!\tau_3}
\!-\!\lr{2}_{\!d_1+1}^{\!\tau_3}\lr{2}_{\!d_2}^{\!\tau_3}\Big)\!\!
\Bigg)$$
for all $d\!\ge\!2$.
The three recursions above and the real divisor relation~\eref{Rdivrel_e} 
determine all numbers~$\lr{\m}_d^{\tau_3}$  with $d\!\ge\!2$.
The low-degree numbers obtained from these recursions are as in Table~4.2.2 in~\cite{Adam}.
Along with~\eref{SNXphi_e}, the numbers~$\lr{\m}_d^{\tau_3}$ determine 
the invariants~\eref{RGWdfn_e} of~$(\P^3,\tau_3)$ 
with all possible insertions~$\mu_i$ in~$H^6(\P^3)$ and~$H^4(\P^3,\R\P^3)$.  
Since only two elements of $H_2(\P^3\!-\!\R\P^3;\Z)$ can be represented by 
a linearly embedded $\P^1\!\subset\!\P^3$,
the invariants~\eref{RGWdfn_e} provide lower bounds for the counts of real curves
in~$(\P^3,\tau_3)$ through line and point constraints;
see Table~1 in~\cite{RealWDVV3}.

\section{The sixfold $\P^1\!\times\!\P^1\!\times\!\P^1$}
\label{P1P1P1_sec}

Let $\tau,\tau_{1,1},\tau_{1,1}'$ be the conjugations of Section~\ref{P1P1_sec}.
The fixed loci of the conjugations 
$$\phi_3\!\equiv\!\tau_{1,1}\!\times\!\tau,
\phi_3'\!\equiv\!\tau_{1,1}'\!\times\!\tau\!:(\P^1)^3\lra(\P^1)^3$$
are the three-torus $(\R\P^1)^3\!\approx\!(S^1)^3$  and
$$\big(\P^1\!\times\!\P^1\big)^{\tau_{1,1}'}\!\times\!(\P^1)^{\tau}
\approx \P^1\!\times\!S^1,$$
respectively.
The groups $H_2^{\phi_3}((\P^1)^3)\!\approx\!H_2((\P^1)^3)$ 
are identified with~$\Z^3$ via the standard generators
$$L_1\equiv\big[\P^1\!\times\!\pt\!\times\!\pt\big], \quad
L_2\equiv\big[\pt\!\times\!\P^1\!\times\!\pt\big], \quad\hbox{and}\quad
L_3\equiv\big[\pt\!\times\!\pt\!\times\!\P^1\big].$$
The group $H_2^{\phi_3'}((\P^1)^3)$ is identified with~$\Z^2$ via the generators
$L_{1,2}\!\equiv\!L_1\!+\!L_2$ and~$L_3$.\\ 

An averager~$\psi$ for both conjugations is given~by 
$$\psi\!:(\P^1)^3\lra(\P^1)^3, \quad 
\psi\big([X_1,Y_1],[X_2,Y_2],[X_3,Y_3]\big)=\big([X_1,Y_1],[X_2,Y_2],[X_3,-Y_3]\big).$$
Each of the two orientations on~$(S^1)^3$ and $\P^1\!\times\!S^1$ is 
compatible with eight Spin-structures and two Spin-structures, respectively.
Both fixed loci have a canonical OSpin-structure~$\os_0$ obtained by trivializing
each factor of~$TS^1$ separately.
It is immediate~that
\BE{P1P1P1_e3}\begin{aligned}
\psi^*\os_0=\ov{\os_0}.
\end{aligned}\EE

\vspace{.15in}

We denote by $e_1,e_2,e_3\!\in\!\Z^3$ the standard basis elements,
by $\1\!\in\!\Z^3$ their sum,
and by $H_1,H_2,H_3\!\in\!H^2((\P^1)^3)$ the Poincare duals of the complex submanifolds
$$\pt\!\times\!\P^1\!\times\!\P^1,\P^1\!\times\!\pt\!\times\!\P^1,
\P^1\!\times\!\P^1\!\times\!\pt\subset(\P^1)^3,$$
respectively.
For an element $m\!=\!(a,b,c)$ of~$\Z^3$, define 
$$|m|=a\!+\!b\!+\!c, \quad
H^m=H_1^aH_2^bH_3^c\in H^{2*}\big((\P^1)^3\big),
\quad mL=aL_1\!+\!bL_2\!+\!cL_3\in H_2\big(\!(\P^1)^3;\Z\big).$$
For an element $d\!\equiv\!(d_1,d_2,d_3)$ of~$(\Z^{\ge0})^3$, 
a tuple $\m\!\equiv\!(m_1,\ldots,m_l)$ of elements of~$(\Z^{\ge0})^3$,
and $\phi\!=\!\phi_3,\phi_3'$, let
$$k_d(\m)= |d|\!+\!l\!-\!\sum_{i=1}^l|m_i|,  \quad
\blr{\m}_{\!d}=\blr{H^{m_1},\ldots,H^{m_l}}_{\!dL}^{\!(\P^1)^3}, \quad
\blr{\m}_{\!d}^{\!\phi}=\blr{H^{m_1},\ldots,H^{m_l}}_{\!dL;\psi}^{\!\phi,\os_0}\,.$$
If in addition $I\!\subset\![l]$, 
let $\m_I$ be the $|I|$-tuple consisting of the entries of~$\m$ indexed by~$I$.
For \hbox{$a_1,a_2,a_3,b\!\in\!\Z^{\ge0}$}, let 
\begin{gather*}
\m_{a_1a_2a_3b}=\big(\underset{a_1}{\underbrace{\1\!-\!e_1,\ldots,\1\!-\!e_1}},
\underset{a_2}{\underbrace{\1\!-\!e_2,\ldots,\1\!-\!e_2}},
\underset{a_3}{\underbrace{\1\!-\!e_3,\ldots,\1\!-\!e_3}},
\underset{b}{\underbrace{\1,\ldots,\1}}\big),\\
N_{d;a_1a_2a_3b}=\blr{\m_{a_1a_2a_3b}}_{\!d}, \quad
N_{d;a_1a_2a_3b}^{\phi}=\blr{\m_{a_1a_2a_3b}}_{\!d}^{\phi}\,.
\end{gather*}
  
\vspace{.15in}

The complex GW-invariants $N_{d;a_1a_2a_3b}$ of $(\P^1)^3$ are preserved by simultaneous 
permutations of the components of the triples~$d$ and~$(a_1,a_2,a_3)$.
Since a degree $d_1L_1\!+\!d_2L_2$ curve is contained in $\P^1\!\times\!\P^1\!\times\!\pt$
for some $\pt\!\in\!\P^1$,
\BE{P1P1P1_e6a} N_{(d_1,d_2,0);a_1a_2a_3b}=\begin{cases}
0,&\hbox{if}~a_1\!+\!a_2\!+\!b\!\neq\!1;\\
N_{d_1,d_2}^{\P^1\times\P^1},&\hbox{if}~(a_1,a_2,a_3,b)\!=\!(0,0,2(d_1\!+\!d_2)\!-\!2,1);\\
d_2N_{d_1,d_2}^{\P^1\times\P^1},&\hbox{if}~(a_1,a_2,a_3,b)\!=\!(1,0,2(d_1\!+\!d_2)\!-\!1,0).
\end{cases}\EE
Since an irreducible degree $d_1L_1\!+\!d_2L_2\!+\!L_3$ curve is the graph of 
a holomorphic map from~$\P^1$ to $\P^1\!\times\!\P^1$,
\BE{P1P1P1_e6b} N_{(d_1,d_2,1);a_1a_2a_3b}=\begin{cases}
0,&\hbox{if}~a_1\!+\!a_2\!+\!b\!<\!3;\\
d_1^{a_2}d_2^{a_1}N_{d_1,d_2}^{\P^1\times\P^1},
&\hbox{if}~a_1\!+\!a_2\!+\!b\!=\!3,\,a_3\!+\!b\!=\!2(d_1\!+\!d_2)\!-\!1.
\end{cases}\EE

\vspace{.15in}

By~\eref{P1P1P1_e3} and Theorem~\ref{WelWDVV3a_thm},
\BE{P1P1P1_e9}\blr{\m}_{\!d}^{\!\phi}=0 \qquad\forall~\phi\!=\!\phi_3,\phi_3',
~d\!\in\!(\Z^{\ge0})^3,\,\m\!\in\!\big(\!(\Z^{\ge0})^3\big)^l~
\hbox{with}~k_d(\m)\!\in\!2\Z.\EE
This implies that
$$\blr{H^{m_1},\ldots,H^{m_l}}_{\!dL;\psi}^{\!\phi,\os_0}
=\blr{H^{m_1},\ldots,H^{m_l}}_{\!dL;\psi}^{\!\phi,{\ov\os_0}}\,,$$ 
i.e.~these invariants do not depend on
an a priori choice of orientation of~$(\!(\P^1)^3)^{\phi}$.
Furthermore,  the invariants $N_{d;a_1a_2a_3b}^{\phi_3}$ are preserved by simultaneous 
permutations of the components of the triples~$d$ and~$(a_1,a_2,a_3)$.
By the same geometric reasoning as in the complex case,
the absolute values of 
the invariants $N_{d;a_1a_2a_3b}^{\phi}$ of $((\P^1)^3,\phi)$
and the invariants $N_{(d_1,d_2),l}^{\vph}$ of $(\P^1)^2$ with
the conjugation $\vph\!=\!\tau_{1,1},\tau_{1,1}'$ corresponding to
$\phi\!\equiv\!\vph\!\times\!\tau$ satisfy the analogues of~\eref{P1P1P1_e6a}
and~\eref{P1P1P1_e6b}.
We describe the relative signs between these invariants in the next paragraph.\\

By~\cite{Jake}, an OSpin-structure~$\os$ on the fixed locus $(\P^1\!\times\!\P^1)^{\vph}$
determines signed counts $N_{(d_1,d_2),l}^{\vph,\os}$ of 
real rational $J$-holomorphic curves  in $(\P^1\!\times\!\P^1,\vph)$
with $|N_{(d_1,d_2),l}^{\vph,\os}|\!=\!|N_{(d_1,d_2),l}^{\vph}|$.
Similarly to the $(\P^1)^3$ case, there is a natural OSpin-structure~$\os_0$ 
on $(\P^1\!\times\!\P^1)^{\vph}$;
the invariants $N_{(d_1,d_2),l}^{\vph,\os_0}$ do not depend on
an a priori choice of orientation of~$(\P^1\!\times\!\P^1)^{\vph}$.
By Theorem~12.1 and Examples~12.3 and~12.4 in~\cite{SpinPin},
\BE{WelComp_e}
N_{(d_1,d_2),l}^{\tau_{1,1},\os_0}=(-1)^{d_1+d_2+l-1}N_{(d_1,d_2),l}^{\tau_{1,1}} 
\quad\hbox{and}\quad
N_{(d,d),l}^{\tau_{1,1}',\os_0}=(-1)^{d+l-1}N_{d,l}^{\tau_{1,1}'}\,.\EE
On the other hand, it is immediate from the definitions of the invariants 
$N_{(d_1,d_2),l}^{\vph,\os_0}$ and $N_{d;a_1a_2a_3b}^{\phi}$
in~\cite{RealWDVV} and~\cite{RealWDVV3}, respectively, 
the CROrient~5a property in Section~7.2 of~\cite{SpinPin}, 
and Corollary~7.7 in~\cite{SpinPin} that 
these invariants satisfy the exact analogues 
of~\eref{P1P1P1_e6a} and~\eref{P1P1P1_e6b}.
The implications of the last two statements for the two involutions~$\phi$
are stated in the respective Sections~\ref{P1P1P1a_subs} and~\ref{P1P1P1b_subs}.\\

Similarly to the $\P^3$ case,
the standard WDVV recursion for the GW-invariants $\lr{\m}_d$ gives
\begin{equation*}\begin{split}
\blr{m_1,m_2,m_3\!+\!e_r,m_4\ldots,m_l}_{\!d}&=\blr{m_1,m_2\!+\!e_r,m_3,m_4,\ldots,m_l}_{\!d}
\!+\!d_r\blr{m_1\!+\!m_3,m_2,m_4,\ldots,m_l}_{\!d}\\
&+\!\!
\sum_{\begin{subarray}{c}d',d''\in(\Z^{\ge0})^3-\{0\}\\ d'+d''=d\end{subarray}}
\!\!\!\Bigg\{\!\!\sum_{(I,J)\in\cP_{3;2}(l)}\!-\!
\sum_{(I,J)\in\cP_{2;3}(l)}\!\!\Bigg\}
~~d_r''\hspace{-.25in}
\sum_{\begin{subarray}{c}i,j\in(\Z^{\ge0})^3-\{0\}\\ i+j=\1\end{subarray}}
\hspace{-.3in}
\blr{\m_I,i}_{\!d'}\blr{j,\m_J}_{\!d''}
\end{split}\end{equation*}
whenever $|m_1|,|m_2|\!\ge\!2$ and $r\!=\!1,2,3$.
Along with the divisor relation, the above identity with 
$$|m_1|\ge |m_2|\ge |m_3|\!+\!1 \ge |m_4|\ge \ldots\ge |m_l|\ge2$$
recursively determines all complex GW-invariants $\lr{\m}_{\!d}$ of~$(\P^1)^3$.

\subsection{The product involution}
\label{P1P1P1a_subs}

By the sentence immediately after~\eref{WelComp_e} and  
the first identity in~\eref{WelComp_e},
\BE{P1P1P1_e11a}\begin{split}
N_{(d_1,d_2,0);a_1a_2a_3b}^{\phi_3,\os_0}&=\begin{cases}
N_{(d_1,d_2),d_1+d_2-1}^{\tau_{1,1}},
&\hbox{if}~a_1,a_2,b\!=\!0,~a_3\!=\!d_1\!+\!d_2\!-\!1;\\
0,&\hbox{otherwise};
\end{cases}\\
N_{(d_1,d_2,1);a_1a_2a_3b}^{\phi_3,\os_0}&=\begin{cases}
0,&\hbox{if}~d_1\!+\!d_2\!+\!a_1\!+\!a_2\!\le\!a_3\!+\!1;\\
-N_{(d_1,d_2),d_1+d_2-2}^{\tau_{1,1}},
&\hbox{if}~a_1,a_2,b\!=\!0,~a_3\!=\!d_1\!+\!d_2\!-\!2;\\
N_{(d_1,d_2),d_1+d_2-1}^{\tau_{1,1}},
&\hbox{if}~a_1,a_2\!=\!0,~b\!=\!1,~a_3\!=\!d_1\!+\!d_2\!-\!2;\\
d_2N_{(d_1,d_2),d_1+d_2-1}^{\tau_{1,1}},
&\hbox{if}~a_1\!=\!1,~a_2,b\!=\!0,~a_3\!=\!d_1\!+\!d_2\!-\!1.
\end{cases}
\end{split}\EE
Along with~\eref{P1P1_e2}, this gives
\BE{P1P1P1_e12a} 
\lr{}_{e_r}^{\phi_3}=1~~\forall\,r\!=\!1,2,3, \qquad 
\lr{}_{e_1+e_2+e_3}^{\phi_3'}=-1\,.\EE

\vspace{.15in}

For $d\!\in\!(\Z^{\ge0})^3$ and 
$\m\!\in\!(\!(\Z^{\ge0})^3)^l$, let
$$\lr{d}_{\m}=\begin{cases}1,&\hbox{if}~k_d(\m)\!=\!2,\,d\!\in\!(2\Z)^3;\\
0,&\hbox{otherwise}.
\end{cases}$$
For the purposes of applying Theorem~\ref{WelWDVV3_thm} in this case,
we take the basis $\{\mu_i^{\st}\}$ for $H^*((\P^1)^3)^{\phi_3}_-$ to 
be~$\{H^m\}$ with $m\!\in\!\{0,1\}^3$.
Taking $\mu_2\!=\!H_r$ with $r\!=\!1,2,3$ in~\ref{Wel12rec3_it},
we obtain
\begin{enumerate}[label=($\phi_31\alph*$),leftmargin=*]

\item If $r\!=\!1,2,3$, $d\!\in\!(\Z^{\ge0})^3\!-\!\{0\}$, 
$l\!\ge\!1$, and $\m\!\equiv\!(m_1,\ldots,m_l)\!\in\!(\!(\Z^{\ge0})^3\!-\!\{0\})^l$
with $k_d(\m)\!\ge\!2$, then
\begin{equation*}\begin{split}
\hspace{-.4in}&\blr{m_1\!+\!e_r,m_2,\ldots,m_l}_{\!d}^{\!\phi_3}
=-2^{l-2}\lr{d}_{\m}d_r\blr{m_1,m_2,\ldots,m_l,\1}_{\!d/2}\\
&\hspace{.5in}-\!\!\!
\sum_{\begin{subarray}{c}d',d''\in(\Z^{\ge0})^3-\{0\}\\ d''+2d'=d\end{subarray}}
d_r'\!\!\!\!\!\!\!\sum_{(I,J)\in\cP(l)}
\!\!\!\!\!\!2^{|I|-1}\!\!\!\!\!\sum_{\begin{subarray}{c}
i,j\in(\Z^{\ge0})^3-\{0\}\\ i+j=\1\end{subarray}}\hspace{-.3in}
\blr{\m_I,i}_{\!d'}\blr{j,\m_J}_{\!d''}^{\!\phi_3}\\
&\hspace{.5in}+\!\!\!\sum_{\begin{subarray}{c}d',d''\in(\Z^{\ge0})^3-\{0\}\\ 
d'+d''=d\end{subarray}}
\sum_{(I,J)\in\cP(l)}\!\!\!
\Bigg(\!d_r''\!\binom{k_d(\m)\!-\!2}{k_{d'}(\m_I)\!-\!1}\!-\!
d_r'\!\binom{k_d(\m)\!-\!2}{k_{d'}(\m_I)}\!\!\Bigg)
\blr{\m_I}_{\!d'}^{\!\phi_3}\blr{\m_J}_{\!d''}^{\!\phi_3}\,.
\end{split}\end{equation*}

\end{enumerate}
Replacing $d$ by $d\!+\!e_r$ above and taking $l\!=\!1$ and $m_1\!=\!e_r$,
we obtain
\begin{enumerate}[label=($\phi_31\alph*$),leftmargin=*]

\setcounter{enumi}{1}

\item If $r\!=\!1,2,3$ and $d\!\in\!(\Z^{\ge0})^3\!-\!\{0,e_r\}$, then 
$$\big(|d|\!-\!1\!-\!2d_r\big)\lr{}_d^{\phi_3}=
\!\!\sum_{\begin{subarray}{c}d',d''\in(\Z^{\ge0})^3\\
|d'|,|d''|\ge2\\ 
d'+d''=d+e_r\end{subarray}}\!\!\!\!\!
d_r'\Bigg(\!d_r''\!\binom{|d|\!-\!1}{|d'|\!-\!1}\!-\!
d_r'\!\binom{|d|\!-\!1}{|d'|}\!\!\Bigg)
\blr{}_{\!d'}^{\!\phi_3}\blr{}_{\!d''}^{\!\phi_3}\,.$$
\end{enumerate}
Along with the $\phi\!=\!\phi_3$ case of~\eref{P1P1P1_e9} and~\eref{P1P1P1_e12a}, 
the two recursions above determine all real GW-invariants~$\lr{\m}_d^{\phi_3}$;
some of them are shown in Table~\ref{P1P1P1a_tbl}.\\

In light of~\eref{P1P1P1_e9}, changing the sign of all invariants~$\lr{\m}_d^{\phi_3}$
with $k_d(\m)\!\equiv\!3$ mod~4 would not invalidate either of the relations of 
Theorem~\ref{WelWDVV3_thm}.
Thus, the last identity in~\eref{P1P1P1_e12a} is not redundant.
Since only one element of the preimage of each~$L_i$ under the first homomorphism
in~\eref{aveG_e0} can be represented by a holomorphic curve,
the invariants~$\lr{\m}_d^{\phi_3}$ provide lower bounds for the counts of real curves
in~$((\P^1)^3,\phi_3)$ through line and point constraints.

\subsection{The twisted involution}
\label{P1P1P1b_subs}

By the sentence immediately after~\eref{WelComp_e} and  
the second identity in~\eref{WelComp_e},
\BE{P1P1P1_e11b}\begin{split}
N_{(d,d,0);a_1a_2a_3b}^{\phi_3',\os_0}&=\begin{cases}
(-1)^dN_{d,2d-1}^{\tau_{1,1}'},
&\hbox{if}~a_1,a_2,b\!=\!0,~a_3\!=\!2d\!-\!1;\\
0,&\hbox{otherwise};
\end{cases}\\
N_{(d,d,1);a_1a_2a_3b}^{\phi_3',\os_0}&=\begin{cases}
0,&\hbox{if}~2d\!+\!a_1\!+\!a_2\!\le\!a_3\!+\!1;\\
(-1)^{d-1}N_{d;2d-2}^{\tau_{1,1}'},
&\hbox{if}~a_1,a_2,b\!=\!0,~a_3\!=\!2d\!-\!2;\\
(-1)^dN_{d,2d-1}^{\tau_{1,1}'},
&\hbox{if}~a_1,a_2\!=\!0,~b\!=\!1,~a_3\!=\!2d\!-\!2;\\
(-1)^ddN_{d,2d-1}^{\tau_{1,1}'},
&\hbox{if}~a_1\!=\!1,~a_2,b\!=\!0,~a_3\!=\!2d\!-\!1.
\end{cases}
\end{split}\EE
Along with~\eref{P1P1_e2}, this gives
\BE{P1P1P1_e12b} 
\lr{}_{e_3}^{\phi_3'}=1, \qquad
\lr{H_1H_2}_{e_1+e_2}^{\phi_3'}=-1, \qquad
\lr{}_{e_1+e_2+e_3}^{\phi_3'}=1\,.\EE

\vspace{.15in}

Let $e_1,e_2\!\in\!\Z^2$ be the standard basis elements and
$$H_1'=\frac12\big(H_1\!+\!H_2), \qquad H_2'=H_3.$$
Thus, $H_1'^2$, $H_1'H_2'$, and $H_1'^2H_2'$ are the Poincare duals
of $L_3/2$, $(L_1\!+\!L_2)/2$, and $\pt/2$, respectively.
For elements $d\!=\!(a,b)$ and $m\!=\!(r,s)$ of $(\Z^{\ge0})^2$, define
$$|d|=a\!+\!b, \qquad m_1=a, \qquad m_2=b, \qquad
H'^m=H_1'^rH_2'^s.$$
If in addition $\m\!\equiv\!(m_1,\ldots,m_l)$ is an element of $(\!(\Z^{\ge0}\!-\!\{0\})^2)^l$, 
let
\begin{gather*}
k_d(\m)=2a\!+\!b\!+\!l\!-\!\sum_{i=1}^l|m_i|,  \quad
\blr{\m}_{\!d}^{\!\phi_3'}=\blr{H'^{m_1},\ldots,H'^{m_l}}_{\!a(L_1+L_2)+bL_3;\psi}^{\!\phi_3',\os_0},\\
\lr{d}_{\m}=\begin{cases}
1,&\hbox{if}~k_d(\m)\!=\!2;\\
0,&\hbox{otherwise}; \end{cases} \quad
\lr{\m}_d=\begin{cases}
\sum\limits_{\begin{subarray}{c}a'+a''=a\\ a',a''\ge0\end{subarray}}
\!\!\!\!\!\blr{H'^{m_1},\ldots,H'^{m_l}}_{\!a'L_1+a''L_2+(b/2)L_3}^{\!(\P^1)^3},&\hbox{if}~b\!\in\!2\Z;\\
0,&\hbox{otherwise}.\end{cases}
\end{gather*}

\vspace{.15in}

For the purposes of applying Theorem~\ref{WelWDVV3_thm} in this case,
we take the basis $\{\mu_i^{\st}\}$ to 
be~$\{H'^m\}$ with $m\!\in\!\{0,1,2\}\!\times\!\{0,1\}$.
Taking $\mu_2\!=\!H_r'$ with $r\!=\!1,2$ in~\ref{Wel12rec3_it},
we obtain
\begin{enumerate}[label=($\phi_3'1\alph*$),leftmargin=*]

\item\label{phi3tw12_it} If $r\!=\!1,2$, $d\!\in\!(\Z^{\ge0})^2\!-\!\{0\}$, 
$l\!\ge\!1$, and $\m\!\equiv\!(m_1,\ldots,m_l)\!\in\!(\!(\Z^{\ge0})^2\!-\!\{0\})^l$
with $k_d(\m)\!\ge\!2$, then
\begin{equation*}\begin{split}
\hspace{-.4in}&\blr{m_1\!+\!e_r,m_2,\ldots,m_l}_{\!d}^{\!\phi_3'}
=-2^{l-1}\lr{d}_{\m}d_r\blr{m_1,m_2,\ldots,m_l,(2,1)}_{\!d}\\
&\hspace{.5in}-\!\!\!
\sum_{\begin{subarray}{c}d',d''\in(\Z^{\ge0})^2-\{0\}\\ d'+d''=d\end{subarray}}
d_r'\!\!\!\!\!\!\!\sum_{(I,J)\in\cP(l)}
\!\!\!\!2^{|I|-1}\!\!\!\!\!\!\!\!\!\!\!\sum_{\begin{subarray}{c}
i,j\in(\Z^{\ge0})^2-\{0\}\\ i+j=(2,1)\end{subarray}}\hspace{-.3in}
\blr{\m_I,i}_{\!d'}\blr{j,\m_J}_{\!d''}^{\!\phi_3'}\\
&\hspace{.5in}+\!\!\!\sum_{\begin{subarray}{c}d',d''\in(\Z^{\ge0})^2-\{0\}\\ 
d'+d''=d\end{subarray}}
\sum_{(I,J)\in\cP(l)}\!\!\!
\Bigg(\!d_r''\!\binom{k_d(\m)\!-\!2}{k_{d'}(\m_I)\!-\!1}\!-\!
d_r'\!\binom{k_d(\m)\!-\!2}{k_{d'}(\m_I)}\!\!\Bigg)
\blr{\m_I}_{\!d'}^{\!\phi_3'}\blr{\m_J}_{\!d''}^{\!\phi_3'}\,.
\end{split}\end{equation*}

\end{enumerate}
Replacing $d\!\equiv\!(a,b)$ by $d\!+\!e_2$ above and taking $r\!=\!2$, $l\!=\!1$, and $m_1\!=\!e_2$,
we obtain
\begin{enumerate}[label=($\phi_3'1\alph*$),leftmargin=*]

\setcounter{enumi}{1}

\item If $(a,b)\in\!(\Z^{\ge0})^2\!-\!\{0,e_2\}$, then 
$$\big(2a\!-\!1\!-\!b\big)\lr{}_{(a,b)}^{\phi_3'}=
\!\!\sum_{\begin{subarray}{c}
a'+a''=a\\ b'+b''=b+1\\ a',b',a'',b''\ge1\end{subarray}}\!\!\!\!\!
b'\Bigg(\!b''\!\binom{2a\!+\!b\!-\!1}{2a'\!+\!b'\!-\!1}\!-\!
b'\!\binom{2a\!+\!b\!-\!1}{2a'\!+\!b'}\!\!\Bigg)
\blr{}_{\!(a',b')}^{\!\phi_3'}\blr{}_{\!(a'',b'')}^{\!\phi_3'}\,.$$
\end{enumerate}
Replacing instead $d$ by $d\!+\!e_1$ in~\ref{phi3tw12_it} and 
taking $r\!=\!1$, $l\!=\!1$, and $m_1\!=\!2e_1$, we obtain
\begin{enumerate}[label=($\phi_3'1\alph*$),leftmargin=*]

\setcounter{enumi}{2}

\item If $(a,b)\!\in\!(\Z^{\ge0})^2\!-\!\{0\}$, then 
$$\big(b\!-\!1\big)\lr{}_{(a,b)}^{\phi_3'}=
\!\!\sum_{\begin{subarray}{c}
a'+a''=a+1\\ b'+b''=b\\ a',b',a'',b''\ge1 \end{subarray}}\!\!\!
\Bigg(\!a'\!\binom{2a\!+\!b\!-\!1}{2a'\!+\!b'\!-\!1}
\!-\!a''\!\binom{2a\!+\!b\!-\!1}{2a'\!+\!b'\!-\!2}
\!\!\Bigg)
\blr{H_1H_2}_{\!(a',b')}^{\!\phi_3'}\blr{}_{\!(a'',b'')}^{\!\phi_3'}\,.$$
\end{enumerate}
Along with the $\phi\!=\!\phi_3'$ case of~\eref{P1P1P1_e9} and~\eref{P1P1P1_e12a}, 
the three recursions above determine all real GW-invariants~$\lr{\m}_d^{\phi_3'}$;
some of them are shown in Table~\ref{P1P1P1b_tbl}.\\

For the same reasons as for the product involution of Section~\ref{P1P1P1b_subs}, 
the last identity in~\eref{P1P1P1_e12b} is not redundant.
Since only one element of the preimage each of $L_1\!+\!L_2$ and~$L_3$ 
under the first homomorphism
in~\eref{aveG_e0} can be represented by a holomorphic curve,
the invariants~$\lr{\m}_d^{\phi_3'}$ provide lower bounds for the counts of real curves
in~$((\P^1)^3,\phi_3')$ through line and point constraints.

\tbl
\begin{tabular}{||c|c|c|c|c|c|c|c|c||}
\hline\hline
& d=1& d=2& d=3& d=4& d=5& d=6& d=7& d=8\\
\hline
$\C$& 1& 1& 12& 620& 87304& 26312976& 14616808192& 13525751027392\\
\hline
l=0& 1& 1&  8& 240& 18264& 2845440& 792731520& 359935488000\\
\hline
l=1& 1& 1&  6& 144& 9096& 1209600& 293758272& 118173265920\\
\hline
l=2& & 1&  4& 80& 4272& 490368& 104600448& 37486448640\\
\hline
l=3& & &  2& 40& 1872& 188544& 35670576& 11463469056\\
\hline
l=4& & &  0& 16& 744& 67968& 11579712& 3367084032\\
\hline
l=5& & &  & 0& 248& 22400& 3538080& 944056320\\
\hline
l=6& & &  & & 64& 6400& 995904& 249999360\\
\hline
l=7& & &  & & 64& 1536& 248976& 61424640\\
\hline
l=8& & &  & & & 1024& 54272& 13643776\\
\hline
l=9& & &  & & & &11776& 2705408\\
\hline
l=10& & &  & & & &-14336& 499712\\
\hline
l=11& & &  & & & & & -280576\\
\hline\hline
\end{tabular}
\newcaption{The counts $N_d^{\P^2}$ of complex genus~0 degree~$d$ curves in~$\P^2$ 
through $3d\!-\!1$ points (the $\C$~line) and 
Welschinger's invariant counts~$N_{d,l}^{\tau_2}$
of real genus~0 degree~$d$ curves in~$\P^2$ through~$l$ conjugate pairs of points 
and $3d\!-\!1\!-\!2l$ real points.}
\endtbl{WelP2_tbl}

\tbl
\begin{tabular}{||c|c|c|c|c|c|c|c||}
\hline\hline
& d=1& d=2& d=3& d=4& d=5& d=6& d=7\\
\hline
$\C$& 1& 12& 3510& 6508640& 43628131782 & 780252921765888 &  \!\!{\small30814236194426422332}\!\!\\
\hline
l=0& 1 & 6 & 576 & 294336 & 493848576 & 2079965454336 &  18546841177030656\\
\hline
l=1& 1 & 4 & 288 & 116352 & 160966656 & 576148930560 &  4464575005261824\\
\hline
l=2& & 2 & 128 & 42624 & 49582080 & 152559783936 &  1035147394547712\\
\hline
l=3& & 0 & 48 & 14208 & 14303232 & 38410813440 &  230362111475712\\
\hline
l=4& & & 16 & 4224 & 3821568 & 9135710208 &  48998058983424\\
\hline
l=5& & & 16 & 1152 & 938496 & 2039070720 &  9915262009344\\
\hline
l=6& & & & 320 & 215040 & 425871360 &  1901347799040\\
\hline
l=7& & & & -256 & 47872 & 83951616 &  345169133568\\
\hline
l=8& & & & & 10496 & 15949824 &  59646984192\\
\hline
l=9& & & & & 26880 & 2998272 &  9935069184\\
\hline
l=10& & &  & & & 630784 &  1624670208\\
\hline
l=11& & &  & & & -2637824 & 270151680\\
\hline
l=12& & &  & & &  &  42536960\\
\hline
l=13& & &  & & & &  500240384\\
\hline\hline
\end{tabular}
\newcaption{The counts $N_{d,d}^{\P^1\!\times\P^1}$ of complex
genus~0 bidegree~$(d,d)$ curves in~$\P^1\!\times\!\P^1$ 
through $4d\!-\!1$ points (the $\C$~line) and 
Welschinger's invariant counts~$N_{d,l}^{\tau_{1,1}'}$
of real genus~0 bidegree~$(d,d)$ curves in~$\P^1\!\times\!\P^1$ 
through~$l$ conjugate pairs of points 
and $4d\!-\!1\!-\!2l$ real points.}
\endtbl{WelP1P1a_tbl}

\tbl
\begin{tabular}{||c|c|c|c|c|c|c|c|c|c|c||}
\hline\hline
& (2,2)& (2,3)& (3,3)& (2,4)& (3,4)& (4,4)& (2,5)& (3,5)& (4,5)& (5,5)\\
\hline
$\C$& 12& 96& 3510& 640&  87544&  6508640&  3840& 1763415& 
  \!\!{\small348005120}\!\!&  \!\!{\small43628131782}\!\!\\
\hline
l=0&  8 & 48 & 1086 & 256 & 18424 & 819200 & 1280 & 268575 & 28312064 & 2082934630 \\
\hline
l=1&  6 & 32 & 606 & 160 & 9256 & 360896 & 768 & 125855 & 11406848 & 756290790 \\
\hline
l=2& 4 & 20 & 318 & 96 & 4432 & 152192 & 448 & 56831 & 4428160 & 265412198 \\
\hline
l=3&  2 & 12 & 158 & 56 & 2032 & 61568 & 256 & 24831 & 1659264 & 90118886 \\
\hline
l=4& & 8 & 78 & 32 & 904 & 24064 & 144 & 10559 & 602496 & 29678982 \\
\hline
l=5& & & 46 & 16 & 408 & 9280 & 80 & 4415 & 213888 & 9532294 \\
\hline
l=6& & & & & 224 & 3712 & 48 & 1887 & 75776 & 3020358 \\
\hline
l=7& & & & & & 1536 &  & 991 & 28160 & 965958 \\
\hline
l=8& & & & & & & & & 13056 & 327974 \\
\hline
l=9& & & & & & & & & & 142758 \\
\hline\hline
\end{tabular}
\newcaption{The counts $N_{a,b}^{\P^1\!\times\P^1}$ of complex genus~0 bidegree~$(a,b)$ 
curves in~$\P^1\!\times\!\P^1$  through $2(a\!+\!b)\!-\!1$ points (the $\C$~line) and 
Welschinger's invariant counts $N_{(a,b),l}^{\tau_{1,1}}$
of real genus~0 bidegree~$(a,b)$ curves in~$\P^1\!\times\!\P^1$ 
through~$l$ conjugate pairs of points 
and $2(a\!+\!b)\!-\!1\!-\!2l$ real points.}
\endtbl{WelP1P1b_tbl}

\tbl
\begin{tabular}{||c|c|c|c|c|c|c|c|c|c|c|c||}
\hline\hline
$d,a$ & 3,2& 4,2& 4,3& 5,2& 5,3& 5,4& 6,2& 6,3& 6,4& 6,5& 7,2\\
\hline
$\C$&  1 & 96 & 1 & 18132 & 640 & 1& 6506400& 401172& 3840& 1& 4059366000\\
\hline
 \text{l=0} & 1 & 48 & 1 & 4584 & 256 & 1& 817920& 71360& 1280& 1& 249486624\\
 \hline
 \text{l=1} & 1 & 32 & 1 & 2412 & 160 & 1& 359616& 34512& 768& 1& 94578912\\
  \hline
 \text{l=2} & 1 & 20 & 1 & 1200 & 96 & 1& 150912& 16000& 448& 1& 34464936\\
  \hline
 \text{l=3} & 1 & 12 & 1 & 564 & 56 & 1& 60288& 7136& 256& 1& 12045432\\
  \hline
 \text{l=4} && 8 & 1 & 248 & 32 & 1& 22784& 3072& 144& 1& 4020816\\
  \hline
 \text{l=5} &&&& 92 & 16 & 1& 8000& 1264& 80& 1& 1271088\\
  \hline
 \text{l=6} &&&& 0 &&& 2432& 448& 48& 1& 373464\\
  \hline
 \text{l=7} &&&&&& & 256&0& && 97352\\
\hline
 \text{l=8} &&&&&& &&&&& 21248\\
\hline
 \text{l=9} &&&&&& &&&&& 13056\\
\hline\hline
\end{tabular}
\newcaption{The  counts $N_{(d,(a))}^1$ of complex genus~0 
degree $dL\!-\!aE$ curves in $\P^2_1$ 
through $3d\!-\!1\!-\!a$ points (the $\C$~line) and 
Welschinger's invariants $N_{(d,(a),()),l}^{1,0}$
of real genus~0 degree $dL\!-\!aE$ curves in $\P^2_{1,0}$ 
through~$l$ conjugate pairs of points and $3d\!-\!1\!-\!a\!-\!2l$ real points.}
\endtbl{WelP2bl1R_tbl}

\tbl
\begin{tabular}{||c|c|c|c|c|c|c|c|c|c||}
\hline\hline
$d,\a$& $4,(\2^2)$& $5,(\2^2)$& $5,(3,2)$& $6,(\2^2)$& $6,(3,2)$& $6,(\3^2)$& $6,(4,2)$
& $7,(\2^2)$\\
\hline
$\C$&   12 & 3510 & 96 & 1558272 & 87544  & 3510& 640& 1108152240\\
\hline
 \text{l=0} &  8 & 1086 & 48 & 229152 & 18424 &1086& 256& 77453856\\
 \hline
 \text{l=1} &  6 & 606 & 32 & 104352 & 9256  & 606& 160& 30056988\\
 \hline
 \text{l=2} &  4 & 318 & 20 & 45312 & 4432  & 318& 96& 11209752\\
 \hline
 \text{l=3} & 2 & 158 & 12 & 18752 & 2032  & 158& 56& 4012308\\
 \hline
 \text{l=4} & & 78 & 8 & 7392 & 904 & 78 & 32& 1374864\\
 \hline
 \text{l=5} & & 46 & & 2784 & 408  & 46& 16& 448812\\
 \hline
 \text{l=6} & & & & 1088 & 224 & && 138056\\
\hline
 \text{l=7} & & & & & & && 38052\\
\hline
 \text{l=8} & & & & & & && 4096\\
\hline\hline
\end{tabular}
\newcaption{The counts $N_{(d,\a)}^2$ of complex genus~0  degree $B_{(d,\a)}$
curves in $\P^2_2$  through $3d\!-\!1\!-\!|\a|$ points (the $\C$~line) and 
Welschinger's invariants $N_{(d,\a,()),l}^{2,0}$
of real genus~0 degree $B_{(d,\a,())}$ curves in $\P^2_{2,0}$ 
through~$l$ conjugate pairs of points and $3d\!-\!1\!-\!|\a|\!-\!2l$ real points.}
\endtbl{WelP2bl2R_tbl}

\tbl
\begin{tabular}{||c|c|c|c|c|c|c|c|c||}
\hline\hline
$d,b$&  4,2& 5,2& 6,2& 6,3& 7,2& 7,3& 8,2& 8,3\\
\hline
$\C$& 12 & 3510 &  1558272 &3510  & 1108152240& 6508640& 1219053648960& 12330654896\\
\hline
 \text{l=0} &  6 & 576 &  88992 &  576& 22823424& 294336& 9282332160& 166440960\\
 \hline
 \text{l=1} &  4 & 288 &  36864 & 288 & 8162688& 116352& 2933701632& 55692288\\
 \hline
 \text{l=2} & 2 & 128 &  14208 & 128 & 2769408& 42624& 888970752& 17639424\\
 \hline
 \text{l=3} & 0 & 48 &  4992 & 48 & 882432& 14208& 256790016& 5240832\\
 \hline
 \text{l=4} && 16 &  1536 & 16 & 259200& 4224& 70027008& 1441536\\
 \hline
 \text{l=5} && 16 &  416 & 16 & 68128& 1152& 17742336& 361728\\
 \hline
 \text{l=6} &&&  288 & & 15936& 320& 4084992& 83584\\
 \hline
 \text{l=7} &&& & & 3616& -256& 848384& 18688\\
 \hline
 \text{l=8} &&&& & -4096& &163840& -8192\\
\hline
 \text{l=9} &&&& & & &-86016&\\
\hline\hline
\end{tabular}
\newcaption{The counts $N_{(d,(b,b))}^2$ of complex genus~0 
degree $dL\!-\!b(E_1\!+\!E_2)$ curves in $\P^2_2$ 
through $3d\!-\!1\!-\!2b$ points (the $\C$~line) and 
Welschinger's invariant counts $N_{(d,(),(b)),l}^{0,1}$
of real genus~0 degree $dL\!-\!bE_1^{\C}$ curves in $\P^2_{0,1}$ 
through~$l$ conjugate pairs of points and $3d\!-\!1\!-\!2b\!-\!2l$ real points.}
\endtbl{WelP2bl1C_tbl}

\tbl
\begin{tabular}{||c|c|c|c|c|c|c|c|c|c|c||}
\hline\hline
\!\!$d,\a$\!\!&    \!\!$5,(\2^3)$\!\!&  \!\!$5,(3,\2^2)$\!\!& 
\!\!$6,(\2^3)$\!\!&  \!\!$6,(3,\2^2)$\!\!& \!\!$6,(\3^2,2)$\!\!
& \!\!$6,(\3^3)$\!\! & \!\!$6,(4,\2^2)$\!\!  & \!\!$7,(\2^3)$\!\!
& \!\!$7,(3,\2^2)$\!\! & \!\!$7,(\3^2,2)$\!\!\\
\hline
$\C$&  620& 12& \!359640\!& \!18132\!& 620& 12& 96& \!296849546\!& \!23133696\!& \!1558272\!\\
\hline
\text{l=0} & 240 & 8 & 62400 & 4584 & 240 & 8 & 48 & 23698434 & 2481632 & 229152 \\
\hline
 \text{l=1} & 144 & 6 & 29520 & 2412 & 144 & 6 & 32 & 9423618 & 1052448 & 104352 \\
\hline
 \text{l=2} & 80 & 4 & 13280 & 1200 & 80 & 4 & 20 & 3598722 & 428032 & 45312 \\
\hline
 \text{l=3} & 40 & 2 & 5680 & 564 & 40 & 2 & 12 & 1318722 & 167040 & 18752 \\
\hline
 \text{l=4} & 16 && 2304 & 248 & 16 & 0 & 8 & 463026 & 62624 & 7392 \\
\hline
 \text{l=5} &&& 848 & 92 & & & & 155378 & 22624 & 2784 \\
\hline
 \text{l=6} &&&&& &&& 50002 & 8128 & 1088 \\
\hline
 \text{l=7} &&&&& &&& 16978 & & \\
\hline\hline
\end{tabular}
\newcaption{The counts $N_{(d,\a)}^3$ of complex genus~0  degree $B_{(d,\a)}$
curves in $\P^2_3$  through $\ell(d,\a)$ points (the $\C$~line) and 
Welschinger's invariants $N_{(d,\a,()),l}^{3,0}$
of real genus~0 degree $B_{(d,\a,())}$ curves in $\P^2_{3,0}$ 
through~$l$ conjugate pairs of points and $\ell(d,\a)\!-\!2l$ real points.}
\endtbl{WelP2bl3R_tbl}

\tbl
\begin{tabular}{||c|c|c|c|c|c|c|c|c|c||}
\hline\hline
\!\!$d,\a\b$\!\!&    \!\!$4,(2),(2)$\!\!&  \!\!$5,(2),(2)$\!\!& 
\!\!$5,(3),(2)$\!\!&  \!\!$6,(2),(2)$\!\!& \!\!$6,(3),(2)$\!\!
& \!\!$6,(2),(3)$\!\! & \!\!$6,(3),(3)$\!\!  & \!\!$6,(4),(2)$\!\!
& \!\!$7,(2),(2)$\!\!\\
\hline
$\C$&  1& 620& 12& 359640& 18132& 620& 12& 96& 296849546\\
\hline
\text{l=0} & 1 & 144 & 8 & 26064 & 2412 & 144 & 8 & 48 & 7330368 \\
\hline
 \text{l=1} & 1 & 80 & 6 & 11328 & 1200 & 80 & 6 & 32 & 2696640 \\
\hline
 \text{l=2} & 1 & 40 & 4 & 4608 & 564 & 40 & 4 & 20 & 943488 \\
\hline
 \text{l=3} && 16 & 2 & 1728 & 248 & 16 & 2 & 12 & 311616 \\
\hline
 \text{l=4} && 0 && 560 & 92 & 0 & 0 & 8 & 95536 \\
\hline
 \text{l=5} &&&& 64 &0&  &  & & 26096 \\
\hline
 \text{l=6} &&&&&&& & & 6160 \\
\hline
 \text{l=7} &&&&&&& & & 3856 \\
\hline\hline
\end{tabular}
\newcaption{The counts $N_{(d,\a\b\b)}^3$ of complex genus~0  degree $B_{(d,\a\b\b)}$
curves in $\P^2_3$  through $\ell(d,\a\b\b)$ points (the $\C$~line) and 
Welschinger's invariants $N_{(d,\a,\b),l}^{1,1}$
of real genus~0 degree $B_{(d,\a,\b)}$ curves in $\P^2_{1,1}$ 
through~$l$ conjugate pairs of points and $\ell(d,\a\b\b)\!-\!2l$ real points.}
\endtbl{WelP2bl1R1C_tbl}

\tbl
\begin{tabular}{||c|c|c|c|c|c|c|c|c|c|c||}
\hline\hline
\!\!$d,\a$\!\!&    \!\!$5,(\2^4)$\!\!&  \!\!$6,(\2^4)$\!\!& 
\!\!$6,(3,\2^3)$\!\!&  \!\!$6,(\3^2,\2^2)$\!\!& \!\!$6,(4,\2^3)$\!\!
& \!\!$7,(\2^4)$\!\! & \!\!$7,(3,\2^3)$\!\!  & \!\!$7,(\3^2,\2^2)$\!\!
& \!\!$7,(\3^3,2)$\!\! & \!\!$7,(\3^4)$\!\!\\
\hline
$\C$&  96& 79416& 3510& 96& 12& 77866800& 5739856& 359640& 18132& 640\\
\hline
\text{l=0} & 48 & 16440 & 1086 & 48 & 8 & 7137408 & 714592 & 62400 & 4584 & 256 \\
\hline
 \text{l=1} & 32 & 8120 & 606 & 32 & 6 & 2912448 & 312208 & 29520 & 2412 & 160 \\
\hline
 \text{l=2} & 20 & 3800 & 318 & 20 & 4 & 1140000 & 130496 & 13280 & 1200 & 96 \\
\hline
 \text{l=3} & 12 & 1688 & 158 & 12 & 2 & 427848 & 52208 & 5680 & 564 & 56 \\
\hline
 \text{l=4} && 728 & 78 & & & 153824 & 20000 & 2304 & 248 & 32 \\
\hline
 \text{l=5} &&& & & & 52688 & 7248 & 848 & & \\
\hline
 \text{l=6} &&& & & & 16512 & & & & \\
\hline\hline
\end{tabular}
\newcaption{The counts $N_{(d,\a)}^4$ of complex genus~0  degree $B_{(d,\a)}$
curves in $\P^2_4$  through $\ell(d,\a)$ points (the $\C$~line) and 
Welschinger's invariants $N_{(d,\a,()),l}^{4,0}$
of real genus~0 degree $B_{(d,\a)}$ curves in $\P^2_{4,0}$ 
through~$l$ conjugate pairs of points and $\ell(d,\a)\!-\!2l$ real points.}
\endtbl{WelP2bl4R_tbl}

\tbl
\begin{tabular}{||c|c|c|c|c|c|c|c|c||}
\hline\hline
\!\!$d,\a,\b$\!\!&    \!\!$5,(\2^2),(2)$\!\!&  \!\!$5,(3,2),(2)$\!\!& 
\!\!$6,(\2^2),(2)$\!\!&  \!\!$6,(3,2),(2)$\!\!& \!\!$6,(\3^2),(2)$\!\!
& \!\!$6,(\2^2),(3)$\!\! & \!\!$6,(3,2),(3)$\!\!& \!\!$7,(\2^2),(2)$\!\!\\
\hline
$\C$&  96&   1& 79416& 3510& 96& 96& 1& 77866800\\
\hline
\text{l=0} & 32 & 1 & 7368 & 606 & 32 & 32 & 1& 2316864 \\
\hline
 \text{l=1} & 20 & 1 & 3360 & 318 & 20 & 20 & 1& 876576 \\
\hline
 \text{l=2} & 12 & 1 & 1440 & 158 & 12 & 12 & 1& 315936 \\
\hline
\text{l=3} & 8 & & 584 & 78 & 8 & 8 & 1& 108040 \\
\hline
\text{l=4} & & & 248 & 46 & & & & 34720\\
\hline
\text{l=5} & & & &  & & & & 9968\\
\hline
\text{l=6} & & & &  & & & & 1152\\
\hline\hline
\end{tabular}
\newcaption{The counts $N_{(d,\a\b\b)}^4$ of complex genus~0  degree $B_{(d,\a\b\b)}$
curves in $\P^2_4$  through $\ell(d,\a\b\b)$ points (the $\C$~line) and 
Welschinger's invariants $N_{(d,\a,\b),l}^{2,1}$
of real genus~0 degree $B_{(d,\a,\b)}$ curves in $\P^2_{2,1}$ 
through~$l$ conjugate pairs of points and $\ell(d,\a\b\b)\!-\!2l$ real points.}
\endtbl{WelP2bl2R1C_tbl}

\tbl
\begin{tabular}{||c|c|c|c|c|c|c|c|c|c||}
\hline\hline
\!\!$d,\b$\!\!&    \!\!$5,(\2^2)$\!\!&  \!\!$6,(\2^2)$\!\!& 
\!\!$6,(3,2)$\!\!&  \!\!$7,(\2^2)$\!\!& \!\!$7,(3,2)$\!\!
& \!\!$7,(\3^2)$\!\! & \!\!$8,(\2^2)$\!\! & \!\!$8,(3,2)$\!\! & \!\!$8,(\3^2)$\!\!\\
\hline
$\C$&  96&   79416& 96& 77866800& 359640& 640& 105128477280& 939726048& 5739856\\
\hline
 \!\!\text{l=0}\!\!&  16&  2640& 16& 625728& 8352& 64& 227372544& 4232448& 65920\\
 \hline
 \!\!\text{l=1}\!\!&  8&   1024& 8& 210240& 3072& 32& 67970688& 1322112& 22528\\
 \hline
 \!\!\text{l=2}\!\!& 4& 352&  4& 65184& 1024& 16& 19175808& 383424& 7168\\
 \hline
 \!\!\text{l=3}\!\!& 4& 112&  4& 18248& 320& 8& 5035008& 102208& 2176\\
 \hline
 \!\!\text{l=4}\!\!&& 80&  & 4640& 96& 0& 1208320& 25344& 640\\
\hline
 \!\!\text{l=5}\!\!&& &  & 1104& -64& &264192& 5984& -128\\
\hline
 \!\!\text{l=6}\!\!&& &  & -1152&&& 53376&  -2304&\\
\hline
 \!\!\text{l=7}\!\!&& &  & &&& -25984&&\\
\hline\hline
\end{tabular}
\newcaption{The counts $N_{(d,\b\b)}^4$ of complex genus~0  degree $B_{(d,\b\b)}$
curves in $\P^2_4$  through $\ell(d,\b\b)$ points (the $\C$~line) and 
Welschinger's invariants $N_{(d,(),\b),l}^{0,2}$
of real genus~0 degree $B_{(d,\b\b)}$ curves in $\P^2_{0,2}$ 
through~$l$ conjugate pairs of points and $\ell(d,\b\b)\!-\!2l$ real points.}
\endtbl{WelP2bl2C_tbl}

\tbl
\begin{tabular}{||c|c|c|c|c|c|c|c|c|c||}
\hline\hline
\!\!$d,\a$\!\!&    \!\!$5,(\2^5)$\!\!&  \!\!$6,(\2^5)$\!\!& 
\!\!$6,(3,\2^4)$\!\!&  \!\!$6,(\3^2,\2^3)$\!\!& \!\!$7,(\2^5)$\!\!
&\!\!$7,(3,\2^4)$\!\! &\!\!$7,(\3^2,\2^3)$\!\! &\!\!$7,(\3^3,\2^2)$\!\!
&\!\!$7,(\3^4,2)$\!\!\\
\hline
$\C$&  12&   16608& 620& 12& 19948176& 1380648& 79416& 3510& 96\\
\hline
 \!\!\text{l=0}\!\!&  8& 4160& 240& 8& 2112480& 201192& 16440& 1086& 48\\
 \hline
 \!\!\text{l=1}\!\!&  6&   2160& 144& 6& 886224& 90856& 8120& 606& 32\\
 \hline
 \!\!\text{l=2}\!\!& 4& 1056&  80& 4& 356076& 39144& 3800& 318& 20\\
 \hline
 \!\!\text{l=3}\!\!&& 480&  40& & 137012& 16104& 1688& 158& 12\\
 \hline
 \!\!\text{l=4}\!\!&&&  && 50568& 6376& 728&&\\
\hline
 \!\!\text{l=5}\!\!&&&  && 18088&&&& \\
\hline\hline
\end{tabular}
\newcaption{The counts $N_{(d,\a)}^5$ of complex genus~0  degree $B_{(d,\a)}$
curves in $\P^2_5$  through $\ell(d,\a)$ points (the $\C$~line) and 
Welschinger's invariants $N_{(d,\a,()),l}^{5,0}$
of real genus~0 degree $B_{(d,\a)}$ curves in $\P^2_{5,0}$ 
through~$l$ conjugate pairs of points and $\ell(d,\a)\!-\!2l$ real points.}
\endtbl{WelP2bl5R_tbl}

\tbl
\begin{tabular}{||c|c|c|c|c|c|c||}
\hline\hline
\!\!$d,\a,\b$\!\!&    \!\!$5,(\2^3),(2)$\!\!&  \!\!$6,(\2^3),(2)$\!\!& 
\!\!$6,(3,\2^2),(2)$\!\!&  \!\!$6,(\3^2,2),(2)$\!\!& \!\!$6,(\2^3),(3)$\!\!
& \!\!$7,(\2^3),(2)$\!\!\\
\hline
$\C$&  12&   16608& 620& 12& 12& 19948176\\
\hline
 \!\!\text{l=0}\!\!&  6&  2004& 144& 6& 6& 720144\\
 \hline
 \!\!\text{l=1}\!\!&  4&   960& 80& 4& 4& 280320\\
 \hline
 \!\!\text{l=2}\!\!& 2& 428&  40& 2& 2& 103948\\
 \hline
 \!\!\text{l=3}\!\!&& 168&  16& & & 36660\\
 \hline
 \!\!\text{l=4}\!\!&&&  &&& 12376\\
\hline
 \!\!\text{l=5}\!\!&&&  &&& 4408\\
\hline\hline
\end{tabular}
\newcaption{The counts $N_{(d,\a\b\b)}^5$ of complex genus~0  degree $B_{(d,\a\b\b)}$
curves in $\P^2_5$  through $\ell(d,\a\b\b)$ points (the $\C$~line) and 
Welschinger's invariants $N_{(d,\a,\b),l}^{3,1}$
of real genus~0 degree $B_{(d,\a,\b)}$ curves in $\P^2_{3,1}$ 
through~$l$ conjugate pairs of points and $\ell(d,\a\b\b)\!-\!2l$ real points.}
\endtbl{WelP2bl3R1C_tbl}

\tbl
\begin{tabular}{||c|c|c|c|c|c|c||}
\hline\hline
\!\!$d,\a,\b$\!\!&    \!\!$5,(2),(\2^2)$\!\!&  \!\!$6,(2),(\2^2)$\!\!& 
\!\!$6,(3),(\2^2)$\!\!&  \!\!$6,(2),(3,\2)$\!\!& \!\!$6,(4),(\2^2)$\!\!
& \!\!$7,(2),(\2^2)$\!\!\\
\hline
$\C$&  12&   16608& 620& 12& 1& 19948176\\
\hline
 \!\!\text{l=0}\!\!&  4&  808& 80& 4& 1& 207744\\
 \hline
 \!\!\text{l=1}\!\!&  2& 336& 40& 2& 1& 72528\\
 \hline
 \!\!\text{l=2}\!\!& 0& 120&  16& 0& 1& 23468\\
 \hline
 \!\!\text{l=3}\!\!&& 16&  0& & & 6804\\
 \hline
 \!\!\text{l=4}\!\!&&&  &&& 1768\\
\hline
 \!\!\text{l=5}\!\!&&&  &&& 1128\\
\hline\hline
\end{tabular}
\newcaption{The counts $N_{(d,\a\b\b)}^5$ of complex genus~0  degree $B_{(d,\a\b\b)}$
curves in $\P^2_5$  through $\ell(d,\a\b\b)$ points (the $\C$~line) and 
Welschinger's invariants $N_{(d,\a,\b),l}^{1,2}$
of real genus~0 degree $B_{(d,\a,\b)}$ curves in $\P^2_{1,2}$ 
through~$l$ conjugate pairs of points and $\ell(d,\a\b\b)\!-\!2l$ real points.}
\endtbl{WelP2bl1R2C_tbl}

\tbl
\begin{tabular}{||c|c|c|c|c|c|c|c|c|c||}
\hline\hline
\!\!$d,\a$\!\!&    \!\!$6,(\2^6)$\!\!&  \!\!$6,(3,\2^5)$\!\!& 
\!\!$7,(\2^6)$\!\!&  \!\!$7,(3,\2^5)$\!\!& \!\!$7,(\3^2,\2^4)$\!\!
& \!\!$7,(\3^3,\2^3)$\!\! & \!\!$7,(\3^4,\2^2)$\!\! & \!\!$7,(4,\2^5)$\!\!
& \!\!$7,(4,3,\2^4)$\!\!\\
\hline
$\C$&  3240&  96& 4974460& 320160& 16608& 620& 12& 3510& 96\\
\hline
 \!\!\text{l=0}\!\!&  1000&  48& 613128& 55168& 4160& 240& 8& 1086& 48\\
 \hline
 \!\!\text{l=1}\!\!&  552&   32& 265074& 25856& 2160& 144& 6& 606& 32\\
 \hline
 \!\!\text{l=2}\!\!& 288& 20&  109532& 11520& 1056& 80& 4& 318& 20\\
 \hline
 \!\!\text{l=3}\!\!&&&  43222& 4864& 480& & & 158&\\
 \hline
 \!\!\text{l=4}\!\!&&&  16240&&& & & &\\
\hline\hline
\end{tabular}
\newcaption{The counts $N_{(d,\a)}^6$ of complex genus~0  degree $B_{(d,\a)}$
curves in $\P^2_6$  through $\ell(d,\a)$ points (the $\C$~line) and 
Welschinger's invariants $N_{(d,\a,()),l}^{6,0}$
of real genus~0 degree $B_{(d,\a,())}$ curves in $\P^2_{6,0}$ 
through~$l$ conjugate pairs of points and $\ell(d,\a)\!-\!2l$ real points.}
\endtbl{WelP2bl6R_tbl}

\tbl
\begin{tabular}{||c|c|c|c|c|c|c|c||}
\hline\hline
\!\!$d,\a,\b$\!\!&    \!\!$6,(\2^4),(2)$\!\!&  \!\!$6,(3,\2^3),(2)$\!\!& 
\!\!$7,(\2^4),(2)$\!\!&  \!\!$7,(3,\2^3),(2)$\!\!& \!\!$7,(\3^2,\2^2),(2)$\!\!
& \!\!$7,(\2^4),(3)$\!\! & \!\!$7,(\3^3,2),(2)$\!\!\\
\hline
$\C$& 3240& 96& 4974460& 320160& 16608& 16608& 620\\
\hline
 \!\!\text{l=0}\!\!&  522& 32& 219912& 22768& 2004& 2004& 144\\
 \hline
 \!\!\text{l=1}\!\!&  266&   20& 88186& 9904& 960& 960& 80\\
 \hline
 \!\!\text{l=2}\!\!& 130& 12&  33644& 4080& 428& 428& 40\\
 \hline
 \!\!\text{l=3}\!\!&&&  12142& 1568& 168& 168& \\
 \hline
 \!\!\text{l=4}\!\!&&&  3984&&&& \\
\hline\hline
\end{tabular}
\newcaption{The counts $N_{(d,\a\b\b)}^6$ of complex genus~0  degree $B_{(d,\a\b\b)}$
curves in $\P^2_6$  through $\ell(d,\a\b\b)$ points (the $\C$~line) and 
Welschinger's invariants $N_{(d,\a,\b),l}^{4,1}$
of real genus~0 degree $B_{(d,\a,\b)}$ curves in $\P^2_{4,1}$ 
through~$l$ conjugate pairs of points and $\ell(d,\a\b\b)\!-\!2l$ real points.}
\endtbl{WelP2bl4R1C_tbl}

\clearpage

\tbl
\begin{tabular}{||c|c|c|c|c|c|c|c||}
\hline\hline
\!\!$d,\a,\b$\!\!&    \!\!$6,(\2^2),(\2^2)$\!\!&  \!\!$6,(3,2),(\2^2)$\!\!& 
\!\!$7,(\2^2),(\2^2)$\!\!&  \!\!$7,(3,2),(\2^2)$\!\!& \!\!$7,(\3^2),(\2^2)$\!\!
& \!\!$7,(\2^2),(3,2)$\!\! & \!\!$7,(3,2),(3,2)$\!\!\\
\hline
$\C$&  3240&   96& 4974460& 320160& 16608& 16608& 620\\
\hline
 \!\!\text{l=0}\!\!&  236&  20& 67608& 8320& 864& 808& 80\\
 \hline
 \!\!\text{l=1}\!\!&  108&   12& 24530& 3344& 392& 336& 40\\
 \hline
 \!\!\text{l=2}\!\!& 52& 8&  8332& 1280& 176& 120& 16\\
 \hline
 \!\!\text{l=3}\!\!&&&  2518& 432& 72& 16& \\
 \hline
 \!\!\text{l=4}\!\!&&&  320&&&& \\
\hline\hline
\end{tabular}
\newcaption{The counts $N_{(d,\a\b\b)}^6$ of complex genus~0  degree $B_{(d,\a\b\b)}$
curves in $\P^2_6$  through $\ell(d,\a\b\b)$ points (the $\C$~line) and 
Welschinger's invariants $N_{(d,\a,\b),l}^{2,2}$
of real genus~0 degree $B_{(d,\a,\b)}$ curves in $\P^2_{2,2}$ 
through~$l$ conjugate pairs of points and $\ell(d,\a\b\b)\!-\!2l$ real points.}
\endtbl{WelP2bl2R2C_tbl}

\tbl
\begin{tabular}{||c|c|c|c|c|c|c||}
\hline\hline
$d,\b$&  $5,(\2^3)$& $6,(\2^3)$&  $6,(3,\2^2)$& $7,(\2^3)$& $7,(3,\2^2)$& $7,(\3^2,2)$\\
\hline
$\C$& 1 & 3240&   1&  4974460& 16608& 12\\
\hline
 \text{l=0}& 1& 78&   1& 15864& 244& 4\\
 \hline
 \text{l=1}& 1& 30&   1& 4794& 88& 2\\
 \hline
 \text{l=2}&& 22&&  1340& 28& 0\\
 \hline
 \text{l=3}&&&&  334& -16& \\
 \hline
 \text{l=4}&&&&  -320&& \\
\hline\hline
\end{tabular}
\newcaption{The counts $N_{(d,\b\b)}^6$ of complex genus~0  degree $B_{(d,\b\b)}$
curves in $\P^2_6$  through $\ell(d,\b\b)$ points (the $\C$~line) and 
Welschinger's invariants $N_{(d,(),\b),l}^{0,3}$
of real genus~0 degree $B_{(d,(),\b)}$ curves in $\P^2_{0,3}$ 
through~$l$ conjugate pairs of points and $\ell(d,\b\b)\!-\!2l$ real points.}
\endtbl{WelP2bl3C_tbl}

\tbl
\begin{tabular}{||c|c|c|c|c|c|c|c|c|c||}
\hline\hline
$(a_1,a_2,a_3)$& (1,0,0)& (1,1,1)& (2,1,0)& (3,0,0)& (2,2,1)& (3,1,1)& (3,2,0)& (4,1,0)& (5,0,0)\\
\hline
$\C$& 12& 110& 48& 0& 672& 192& 192& 0& 0\\
\hline
b=0& 6& -12& -8& 0& 16& 8& 8& 0& 0\\
\hline
b=1& -4& 6& 4& 0&&&&&\\
\hline
b=2& 2& & & &&&&&\\
\hline\hline
\end{tabular}
\newcaption{The counts $N_{(2,2,2);2\a}$, with $\a\!\equiv\!(a_1,a_2,a_3)$ 
of complex genus~0  degree~$(2,2,2)$ curves in~$(\P^1)^3$ through 
$2a_1,2a_2,2a_3$ generic representatives for the standard line classes $L_1,L_2,L_3$ 
and $6\!-\!|\a|$ general points (the $\C$~line)
and the signed count $N_{(2,2,2);\a b}^{\phi_3}$
of real genus~0 degree~$(2,2,2)$ curves in~$(\P^1)^3$ 
through $a_1,a_2,a_3$ conjugate pairs of 
generic representatives for $L_1,L_2,L_3$,
$b$ conjugate pairs of points, and $6\!-\!|\a|\!-\!2b$ real points.}
\endtbl{P1P1P1a_tbl}

\clearpage

\tbl
\begin{tabular}{||c|c|c|c|c|c|c|c|c|c|c|c|c||}
\hline\hline
$(a_1,a_2)$
& (1,0)& (0,1)& (2,1)& (1,2)& (3,0)& (0,3)& (3,2)& (2,3)& (4,1)& (1,4)& (5,0)& (0,5)\\
\hline
$\C$& 16& 12& 110& 48& 64& 0& 788& 192& 672& 0& 256& 0\\
\hline b=0&  4& 4& -4& -4& 0& 0& 0& 0& -4& 0& -4& 0\\
\hline
b=1& -2& -2& 0& 0& -2& 0&&&&&&\\
\hline
b=2& 0& 0&&&&&&&&&&\\
\hline\hline
\end{tabular}
\newcaption{The counts $N_{(2,2,2);(a_1,a_1,2a_2)}$ 
of complex genus~0  degree~$(2,2,2)$ curves in~$(\P^1)^3$ through 
$a_1,a_1,2a_2$ generic representatives for the standard line classes $L_1,L_2,L_3$ 
and $6\!-\!a_1\!-\!a_2$ general points (the $\C$~line)
and the signed count $N_{(2,2,2);a_10a_2b}^{\phi_3'}$
of real genus~0 degree~$(2,2,2)$ curves in~$(\P^1)^3$ 
through $a_1$ conjugate pairs of generic representatives for $L_1,L_2$,
$a_2$~conjugate pairs of generic representatives for~$L_3$,
$b$ conjugate pairs of points, and $6\!-\!a_1\!-\!a_2\!-\!2b$ real points.}
\endtbl{P1P1P1b_tbl}

\vspace{.5in}

\noindent
{\it Department of Mathematics, Stony Brook University, Stony Brook, NY 11794\\
xujia@math.stonybrook.edu, azinger@math.stonybrook.edu}

\end{document}